\documentclass[12pt,francais, reqno]{smfart}

\usepackage{amsmath}
\usepackage{amssymb}
\usepackage{amsthm}
\usepackage{a4wide}
\usepackage[frenchb,english]{babel}
\usepackage{mathrsfs}

\numberwithin{equation}{section}

\newcommand{\dd}{\textup{d}}

\def\al{\alpha}

\def\lam{\lambda}

\def\dd{\textup{d}}

\newcommand{\moins}{\setminus}
\newcommand{\card}{{\rm Card}} 
\newcommand{\Id}{{\rm Id}}
\newcommand{\Z}{\mathbb Z}
\newcommand{\Q}{\mathbb Q}
\newcommand{\R}{\mathbb R}

\newcommand{\gdo}{\mathcal{O}}
\newcommand{\Deltati}{\widetilde{\Delta}}
\newcommand{\indso}[2]{{\tiny {\begin{array}{c} #1 \\ #2 \end{array}}}}

\newcommand{\cjs}{C\bigg[\,\begin{matrix} s_1, \ldots, s_p\\j_1, \ldots, j_p\end{matrix}\,\bigg]}
\newcommand{\cjsi}{C\bigg[\,\begin{matrix} s_1, \ldots, s_p\\j_1, \ldots,  j_{i-1}, n-j_i, j_{i+1}, \ldots, j_p\end{matrix}\,\bigg]}
\newcommand{\cjsgamma}{C\bigg[\,\begin{matrix} s_{\gamma(1)}, \ldots , s_{\gamma(p)} \\ j_{\gamma(1)},\ldots, j_{\gamma(p)} \end{matrix}\,\bigg]}
\newcommand{\cjsun}{C\bigg[\,\begin{matrix} 1, s_2,  \ldots, s_p\\j_1, j_2, \ldots, j_p\end{matrix}\,\bigg]}
\newcommand{\cjsinv}{C\bigg[\,\begin{matrix} s_1, \ldots, s_{i-1}, 1, s_{i+1},  \ldots, s_p\\j_1,  \ldots, j_{i-1}, j_i, j_{i+1}, \ldots,  j_p\end{matrix}\,\bigg]}
\newcommand{\cjsparsix}{C\bigg[\,\begin{matrix} s_1, \ldots, s_{\ell-1}, s_\ell, s_{\ell+1}, \ldots, s_p\\j_1, \ldots,  j_{\ell-1}, n-j_\ell, j_{\ell+1}, \ldots, j_p\end{matrix}\,\bigg]}
\newcommand{\cjsparsixbis}{C\bigg[\,\begin{matrix} s_1, \ldots, s_p\\ \eps_1 \cdot j_1, \ldots, \eps_p \cdot j_p\end{matrix}\,\bigg]}
\newcommand{\cjspr}{C\bigg[\,\begin{matrix} s'_1, \ldots, s'_p\\j'_1, \ldots, j'_p\end{matrix}\,\bigg]}

\newcommand{\orb}{\Omega_{\jsoul, \ssoul}}

\newcommand{\spp}{{\mathfrak S}_p}
\newcommand{\sppmd}{{\mathfrak S}_{p-2}}
\newcommand{\strois}{\mathfrak{S}_3}

\newcommand{\zdzp}{(\Z / 2 \Z)^p}
\newcommand{\zdzpmd}{(\Z / 2 \Z)^{p-2}}
\newcommand{\zdzt}{(\Z / 2 \Z)^3}
\newcommand{\zdzd}{(\Z / 2 \Z)^2}
\newcommand{\zdz}{\Z / 2 \Z}

\newcommand{\psd}{\rtimes}

\newcommand{\qjs}{Q_{\underline j, \underline s}}

\newcommand{\eps}{\varepsilon}
\newcommand{\epsunz}{\eps_1^0}
\newcommand{\epstetaz}{\eps_\teta^0}
\newcommand{\sigmati}{\widetilde \sigma}

\newcommand{\gdoeps}{\gdo_\eps (N^{-1+\eps})}
\newcommand{\epssoul}{\underline{\eps}}
\newcommand{\jsoul}{\underline{j}}
\newcommand{\jprsoul}{\underline{j'}}
\newcommand{\jtisoul}{\underline{j} + (\eps_1, 0, \ldots, 0)}

\newcommand{\ssoul}{\underline{s}}
\newcommand{\sprsoul}{\underline{s'}}

\newcommand{\jprun}{j'_1+\frac{\eps_1-1}{2}}

\newcommand{\zerop}{\{0,\ldots, p\}}
\newcommand{\zeron}{\{0,\ldots, n\}}
\newcommand{\zeronmu}{\{0,\ldots, n-1\}}
\newcommand{\unn}{\{1,\ldots, n\}}
\newcommand{\una}{\{1,\ldots, A\}}
\newcommand{\unp}{\{1,\ldots, p\}}
\newcommand{\unpmu}{\{1,\ldots, p-1\}}
\newcommand{\unpmd}{\{1,\ldots, p-2\}}
\newcommand{\deuxp}{\{2,\ldots, p\}}
\newcommand{\unq}{\{1,\ldots, q\}}
\newcommand{\unqpr}{\{1,\ldots, q'\}}
\newcommand{\unqsec}{\{1,\ldots, q''\}}

\newcommand{\tsig}{t_{\sigma}}
\newcommand{\tsigti}{t_{\sigmati}}
\newcommand{\desj}{\Delta_{\epssoul} ^\sigma (\jsoul)}
\newcommand{\sesj}{S_{\epssoul} ^\sigma (\jsoul)}
\newcommand{\sesjti}{\widetilde{S}_{\epssoul} ^\sigma (\jsoul)}
\newcommand{\sesjtisigti}{\widetilde{S}_{\epssoul} ^{\sigmati} (\jsoul)}
\newcommand{\sesjtirondjtisoul}{\widetilde{S}_{\epssoul} ^{\sigma \circ \transpotsig} (\jtisoul)}
\newcommand{\sesjtijti}{\widetilde{S}_{\epssoul} ^\sigma (\jtisoul)}
\newcommand{\sesigtj}{\widetilde{S}_{\epssoul} ^{\sigt} (\jsoul)}

\newcommand{\setiphimu}{\widetilde{S}_{\epssoul} ^{\Phi^{-1}(t,\teta,\gamma)} (\jsoul)}
\newcommand{\setiphimujti}{\widetilde{S}_{\epssoul} ^{\Phi^{-1}(t,\teta,\gamma)} (\jtisoul)}

\newcommand{\astgj}{{\mathcal A}_{\epssoul} ^{\teta, \gamma} (\jsoul)}
\newcommand{\astgseul}{{\mathcal A}_{\epssoul} ^{\teta, \gamma}  }
\newcommand{\astgjti}{{\mathcal A}_{\epssoul} ^{\teta, \gamma} (\jtisoul)}
\newcommand{\bstgj}{{\mathcal B}_{\epssoul, i} ^{\teta, \gamma} (\jsoul)}
\newcommand{\bstgjti}{{\mathcal B}_{\epssoul, i} ^{\teta, \gamma} (\jtisoul)}
\newcommand{\bstgseul}{{\mathcal B}_{\epssoul, i} ^{\teta, \gamma}  }
\newcommand{\bstgjprime}{{\mathcal B}_{\epssoul, i} ^{\teta', \gamma'} (\jsoul)}

\newcommand{\sigt}{\sigma_t}

\newcommand{\tetasigt}{\teta_{\sigt}}
\newcommand{\tetasig}{\teta_{\sigma}}
\newcommand{\teta}{\vartheta}

\newcommand{\phisig}{\phii_{\sigma}}
\newcommand{\phisigt}{\phii_{\sigt}}
\newcommand{\phii}{\varphi}

\newcommand{\psisig}{\psii_{\sigma}}
\newcommand{\psisigt}{\psii_{\sigt}}
\newcommand{\psiteta}{\psii_{\teta}}
\newcommand{\psitetapr}{\psii_{\teta'}}
\newcommand{\psii}{\psi}

\newcommand{\transpotsig}{ (\tsig-1 \, \, \tsig)} 
\newcommand{\albe}{(\alpha \, \, \ldots \, \, \beta)}
\newcommand{\albelongun}{(\alpha \, \, \alpha+1 \, \,  \ldots \, \, \beta-1 \, \, \beta)}
\newcommand{\albelongde}{(\alpha \, \, \alpha-1 \, \,  \ldots \, \, \beta+1 \, \, \beta)}

\newcommand{\croix}{\times}
\newcommand{\calr}{{\mathscr{R}}}

\newcommand{\N}{\mathbb{N}}
\newcommand{\Netoile}{{\N ^*}}

\newcommand{\zetaetoile}{\zeta _*}
\newcommand{\zetaetoileantisym}{\zeta _* ^{{\rm as}}}
\newcommand{\zetaantisym}{\zeta^{{\rm as}}}  

\newcommand{\odu}{o(1)}
\newcommand{\opdu}{o_p(1)}

\newcommand{\calf}{{\mathscr{F}}}
\newcommand{\cale}{{\mathscr{E}}}
\newcommand{\Card}{{\rm Card}}

\newcommand{\ejsr}{E_{j,s}^{(r)}}
\newcommand{\ejsri}{E_{j_i,s_i}^{(r_i)}}

\title[]{Ph\'enom\`enes de sym\'etrie dans des formes lin\'eaires en polyz\^etas}
\author[]{J. Cresson, S. Fischler et T. Rivoal} 
\date{\today}

\newtheorem{Th}{Th\'eor\`eme}
\newtheorem{Lemme}{Lemme}
\newtheorem{Prop}{Proposition}
\newtheorem{Cor}{Corollaire}
\newtheorem{Conj}{Conjecture}
\newtheorem{remark}{Remarque}[subsection]

\begin{document}
\maketitle

\setcounter{tocdepth}{2}
\baselineskip 6mm

\begin{abstract}
On donne deux g\'en\'eralisations, en profondeur quelconque, du ph\'enom\`ene
de sym\'etrie utilis\'e par Ball-Rivoal pour d\'emontrer qu'une infinit\'e de valeurs de la fonction $\zeta$ de Riemann aux entiers impairs sont irrationnelles. Ces g\'en\'eralisations concernent des s\'eries multiples de type hyperg\'eom\'etrique qui s'\'ecrivent comme   formes lin\'eaires en certains polyz\^etas. La preuve utilise notamment la r\'egularisation des polyz\^etas \`a divergence logarithmique.
\end{abstract}

\begin{altabstract}
We give two generalizations, in arbitrary depth, of the symmetry phenomenon used by  Ball-Rivoal to prove that infinitely many values of Riemann  $\zeta$ function at odd integers are irrational. These generalizations concern multiple series of hypergeometric type, which can be written as linear forms in some specific multiple zeta values. The proof makes use of the regularization procedure for multiple zeta values with logarithmic divergence. 
\end{altabstract}

\tableofcontents

\section{Introduction}\label{sec:1}
Une g\'en\'eralisation 
de la fonction z\^eta de Riemann $\zeta(s)$ est 
donn\'ee par les s\'eries {\em polyz\^etas}, d\'efinies pour tout entier 
$p\ge 1$ 
et tout $p$-uplet $\underline{s}=(s_1, s_2, \dots, s_p)$ 
d'entiers $\ge 1$,
avec $s_1\ge 2$, par
$$
\zeta(s_1, s_2, \ldots, s_p)=
\sum_{k_1> k_2>\ldots > k_p\ge 1}
\frac{1}{k_1^{s_1}k_2^{s_2}\ldots k_p^{s_p}}.
$$
Les entiers $p$ et  
$s_1+s_2+\ldots+s_p$ sont respectivement la profondeur et 
le poids de $\zeta(s_1, s_2, \ldots, s_p)$. 
On voit naturellement appara{\^\i}tre les polyz\^etas lorsque, par exemple,  
on consid{\`e}re les produits des valeurs de la fonction 
z\^eta : 
on a  $\zeta(n)\zeta(m)=\zeta(n+m)+\zeta(n,m)+\zeta(m,n)$, 
ce qui permet en 
quelque sorte  de {\og lin{\'e}ariser\fg} ces produits. 
En dehors de quelques identit\'es telles que $\zeta(2,1)=\zeta(3)$
 (due \`a Euler), 
la nature arithm\'etique de ces s\'eries est aussi peu connue que 
celle des nombres $\zeta(s)$. 
Cependant, l'ensemble des nombres $\zeta(\underline s)$ poss\`ede 
une tr\`es riche 
structure alg\'ebrique assez bien comprise, au moins 
conjecturalement (voir~\cite{MiW}). 
Par exemple, on peut s'int{\'e}resser
 aux $\Q$-sous-espaces vectoriels $\mathcal{Z}_p$ de
 $\mathbb{R}$, engendr{\'e}s par les $2^{p-2}$ polyz{\^e}tas de 
poids $p\ge 2$ :  $\mathcal{Z}_2=\Q\zeta(2)$, 
$\mathcal{Z}_3=\Q\zeta(3)+\Q\zeta(2,1)$, 
$\mathcal{Z}_4=\Q\zeta(4)+\Q\zeta(3,1)+
\Q\zeta(2,2)+
\Q\zeta(2,1,1)$, 
etc. Posons $v_p=\textup{dim}_{\Q}(\mathcal{Z}_p)$. 
On a alors la conjecture suivante, dont le point $(i)$ est d\^u 
\`a Zagier et le point $(ii)$ \`a Goncharov.

\begin{Conj}\label{conj:KZ} 
$(i)$ 
Pour tout entier $p\ge 2$, on a $v_p=c_p$, o{\`u} l'entier 
$c_p$ est d{\'e}fini par 
la r{\'e}currence  lin\'eaire 
$c_{p+3}=c_{p+1}+c_{p}$, 
avec $c_0=1$, $c_1=0$ et $c_2=1$.

 $(ii)$ Les $\Q$-espaces vectoriels $\Q$ et 
$\mathcal{Z}_p$ ($p\ge 2)$, sont en somme
directe.
\end{Conj}
La suite $(v_p)_{p\ge 2}$ devrait donc cro{\^\i}tre comme 
$\al^p$ (o{\`u} $\al\approx
1,3247$ est racine du 
polyn{\^o}me $X^3-X-1$), ce qui est  bien plus petit que 
$2^{p-2}$. Il y a 
donc conjecturalement beaucoup de relations lin{\'e}aires entre les polyz{\^e}tas 
de m{\^e}me poids 
et aucune en poids diff{\'e}rents : 
dans cette direction, un  
th{\'e}or{\`e}me de Goncharov~\cite{Goncharov} et Terasoma~\cite{terasoma}
affirme que l'on a  $v_p\le c_p$ pour tout entier $p\ge 2$. Il
reste donc {\`a} montrer l'in{\'e}galit{\'e} inverse pour montrer $(i)$,  mais aucune 
minoration non
triviale de $v_p$ n'est connue {\`a} ce jour : m\^eme si les relations classiques
donnent  $v_2=v_3=v_4=1$,  
on est bloqu{\'e} d{\`e}s l'{\'e}galit{\'e} $v_5=2$, qui 
est {\'e}quivalente 
{\`a} l'irrationalit{\'e} toujours inconnue de
$\zeta(5)/(\zeta(3)\zeta(2))$. 
Plus g\'en\'eralement, un des int\'er\^ets de la 
Conjecture~\ref{conj:KZ} est d'impliquer la 
suivante. 

\begin{Conj}\label{conj:zeta alg inde} 
Les nombres $\,\pi, \zeta(3), \zeta(5), \zeta(7), \zeta(9),$ etc,
sont alg{\'e}brique\-ment ind{\'e}pen\-dants sur $\Q$.
\end{Conj}

Cette conjecture semble actuellement totalement hors de 
port\'ee. Un certain nombre de r\'esultats diophantiens ont 
n\'eanmoins 
\'et\'e  obtenus en profondeur 1, c'est-\`a-dire dans le cas 
de la 
fonction z\^eta de Riemann (voir \cite{SFBou}) : 
\begin{itemize}
\item[$(i)$] Le nombre $\zeta(3)$ est irrationnel (Ap\'ery~\cite{Apery}) ; 
\item[$(ii)$] La dimension de l'espace vectoriel engendr\'e sur $\Q$ par 1, 
$\zeta(3)$, $\zeta(5), \ldots, \zeta(A)$ (avec $A$ impair) cro{\^\i}t au moins comme $\log(A)$ 
(\cite{BR, RivoalCRAS}) ; 
\item[$(iii)$] Au moins un des quatre nombres $\zeta(5), \zeta(7), \zeta(9), \zeta(11)$ est irrationnel 
(Zudilin~\cite{Zudilinonze}).
\end{itemize}
Ces r\'esultats peuvent \^etre obtenus par l'\'etude de  
certaines s\'eries de la forme 
\begin{equation}\label{eq011}
\sum_{k=1}^{\infty} \frac{P(k)}{(k)_{n+1}^A}
\end{equation}
avec $P(X)\in\Q[X]$, $n\ge0$, $A\ge 1$ ; on utilise ici le symbole de Pochhammer 
d\'efini par $(k)_{\al} = k(k+1)\ldots (k+\al-1)$. Ces s\'eries s'expriment 
 comme combinaisons  lin\'eaires sur $\Q$  de 1 et des valeurs de z\^eta aux 
entiers. Le point crucial est que, dans ces combinaisons lin\'eaires, figurent 
seulement {\em certaines} valeurs de la fonction z\^eta : $\zeta(3)$ dans le cas $(i)$, des 
valeurs $\zeta(s)$ avec $s$ impair dans les cas $(ii)$ et $(iii)$. Ceci provient 
(dans les deux derniers cas, et aussi dans certaines preuves de $(i)$) d'une 
propri\'et\'e de sym\'etrie li\'ee \`a l'aspect (tr\`es) bien \'equilibr\'e\footnote{Dans cet article, nous utilisons
indiff\'eremment les mots {\em (very) well-poised} ou leur traduction fran\c caise {\em (tr\`es) bien \'equilibr\'e}.}
 de 
la s\'erie \eqref{eq011} (voir \cite{BR} ou \cite{RivoalCRAS}) : 

\begin{Th} \label{thwp}
Soit $P \in \Q[X]$ de degr\'e au plus $A(n+1)-2$, tel que 
$$P(-n-X) = (-1)^{A(n+1)+1} P(X).$$
Alors la s\'erie \eqref{eq011} est une combinaison  lin\'eaire, \`a coefficients 
rationnels, de 1 et des valeurs $\zeta(s)$ pour $s$ entier impair compris entre 3 et $A$.
\end{Th}

Le but de cet article est de donner deux g\'en\'eralisations, en profondeur quelconque, 
de ce ph\'enom\`ene de sym\'etrie. Nous esp\'erons que ces g\'en\'eralisations ouvriront 
la porte \`a des r\'esultats diophantiens (d'irrationalit\'e ou d'ind\'ependance lin\'eaire) sur les polyz\^etas qui interviennent (voir \S \ref{subsec23}). 

\bigskip

Notre premier r\'esultat (d\'emontr\'e au paragraphe \ref{sec6}) concerne des 
sommes {\em d\'ecoupl\'ees}, c'est-\`a-dire portant sur tous les $p$-uplets $(k_1, \ldots, k_p) \in \Netoile ^p$ :

\begin{Th} \label{thdecouple}
Soient $p \geq 1$, $n \geq 0$ et $A \geq 1$ des entiers.
Soit $P \in \Q[X_1, \ldots, X_p]$ un polyn\^ome de degr\'e $\leq A(n+1)-2$ par rapport \`a chacune des variables, tel que
\begin{multline*}
\quad P(X_1,\ldots, X_{j-1}, -X_j-n, X_{j+1}, \ldots, X_p ) 
\\
= (-1)^{A(n+1)+1} P( X_1,\ldots, X_{j-1}, X_j, X_{j+1}, \ldots, X_p )\quad 
\end{multline*}
pour tout $j \in \unp$. Alors la somme multiple
\begin{equation} \label{eqdecouple}
\sum_{k_1,  \ldots ,  k_p \ge 1} 
\frac{P(k_1, \ldots, k_p)}{(k_1)_{n+1}^{A} \ldots (k_p)_{n+1}^{A} }
\end{equation}
est un polyn\^ome  \`a coefficients rationnels, de degr\'e au plus $p$,  en les  $\zeta(s)$, pour $s$ entier impair compris entre 3 et $A$.
\end{Th}
Par exemple, lorsque $A=3$ ou $A=4$, cette somme est un polyn\^ome en $\zeta(3)$. Quand on prend 
$p=1$, on retrouve exactement le th\'eor\`eme \ref{thwp} (quel que soit $A$).

La preuve du th\'eor\`eme  \ref{thdecouple} consiste essentiellement (apr\`es avoir d\'ecompos\'e 
la fraction rationnelle en \'el\'ements simples)  \`a s\'eparer la somme multiple en un produit 
de $p$ sommes simples auxquelles on applique le th\'eor\`eme \ref{thwp}. Elle utilise aussi un 
processus de r\'egularisation, dans une situation simple et \'el\'ementaire. 

\bigskip

L'inconv\'enient principal du th\'eor\`eme \ref{thdecouple}, du point de vue des applications 
\'eventuelles, est le fait que la somme sur $k_1$, \ldots, $k_p$ soit d\'ecoupl\'ee. Cet
inconv\'enient est visible par trois aspects que nous d\'ecrivons maintenant.

Tout d'abord, les s\'eries d\'ecoupl\'ees donnent toujours des polyn\^omes en valeurs de $\zeta$ en des entiers, m\^eme quand on omet l'hypoth\`ese de sym\'etrie du  th\'eor\`eme \ref{thdecouple}. Cette remarque, qui d\'ecoule de la preuve du  th\'eor\`eme \ref{thdecouple} (voir \S \ref{sec6}), montre que les polyz\^etas ne peuvent pas intervenir r\'eellement dans ce cadre.

\bigskip

Ensuite, consid\'erons  la s\'erie de Ball  
$$
S_n=n!^2 \sum_{k=1}^{\infty} (k + \frac{n}{2}) \frac{(k-n)_{n}(k+n+1)_n}{(k)_{n+1}^4}.
$$
Pour tout entier $n$, $S_n$ est une forme lin\'eaire en $1$ et $\zeta(3)$ ;  
cela se d\'eduit du th\'eor\`eme \ref{thwp}. Elle co{\"\i}ncide 
exactement avec les formes lin\'eaires  qui ont permis 
\`a Ap\'ery de d\'emontrer l'irrationalit\'e de 
$\zeta(3)$ ; sans rentrer dans les d\'etails, indiquons que 
cette co{\"\i}ncidence n'est pas du tout \'evidente et qu'elle 
est la premi\`ere application de la conjecture 
des d\'enominateurs prouv\'ee dans~\cite{KR}. 
Pour tout entier $p\ge 1$, la s\'erie $S_n^p$ est \'evidemment une s\'erie d\'ecoupl\'ee de la forme consid\'er\'ee dans 
le th\'eor\`eme~\ref{thdecouple} avec 
\begin{multline*}
P(X_1, \ldots, X_p)
\\
= n!^{2p}  (X_1  + \frac{n}{2})\ldots  (X_p + \frac{n}{2}) (X_1-n)_{n}\ldots (X_p-n)_{n} (X_1+n+1)_n \ldots (X_p+n+1)_n 
\end{multline*}
et $A = 4$. Ainsi,  $S_n^p$ est un polyn\^ome
 en $\zeta(3)$ de degr\'e (au plus) $p$, dont on pourrait a priori esp\'erer d\'eduire la transcendance de $\zeta(3)$. Pourtant, 
$S_n^p$ ne contient pas plus d'information diophantienne que $S_n$ et elle ne donne que l'irrationalit\'e de $\zeta(3)$. 

\bigskip

Enfin, les  sommes multiples qui apparaissent dans les preuves d'irrationalit\'e sont plut\^ot de la forme
\begin{equation} \label{eq012}
\sum_{k_1 \geq  \ldots \geq   k_p \ge 1} 
\frac{P(k_1, \ldots, k_p)}{(k_1)_{n+1}^{A} \ldots (k_p)_{n+1}^{A} }, 
\end{equation}
c'est-\`a-dire que la somme porte sur des variables ordonn\'ees ;  c'est \`a ce genre de s\'eries 
que s'applique l'algorithme de \cite{CFRalgo}.  Par exemple, lorsque $p=2$, $A=2$ et
$$
P(X_1, X_2) = n! (X_1-X_2+1)_n (X_2-n)_n(X_2)_{n+1},
$$
Sorokin~\cite{SorokinApery} d\'emontre que la somme \eqref{eq012}  
est exactement~\footnote{Quand on applique l'algorithme 
de \cite{CFRalgo}, on trouve une forme lin\'eaire en 1 et $\zeta(2,1)$ ; il faut alors utiliser 
la relation $\zeta(2,1) = \zeta(3)$. De plus, Sorokin travaille \`a l'aide d'une expression int\'egrale 
alternative de cette somme.} la forme lin\'eaire en 1 et $\zeta(3)$ utilis\'ee par Ap\'ery 
dans sa preuve d'irrationalit\'e. Plus g\'en\'eralement, une conjecture de Vasilyev~\cite{Vasilyev}
affirmait qu'une certaine int\'egrale multiple, \'egale \`a la s\'erie 
\begin{equation}\label{eq:vasilyev}
n!^{p-\varepsilon}\sum_{k_1 \geq  \cdots \geq   k_p \ge 1} 
\frac{(k_1-k_2+1)_n \ldots(k_{p-1}-k_p+1)_n (k_p-n)_n }{(k_1)_{n+1}^{2} \ldots (k_{p-1})_{n+1}^{2}(k_p)_{n+1}^{2-\varepsilon} }, 
\end{equation}
est une forme lin\'eaire rationnelle en les valeurs de z\^eta aux entiers $\ge 2$ 
de la m\^eme parit\'e que $\varepsilon\in\{0,1\}$. La formulation int\'egrale de cette 
conjecture a \'et\'e d\'emontr\'ee dans~\cite{Zudilinservice} et une version 
raffin\'ee dans~\cite{KR} : la m\'ethode consiste \`a prouver que la s\'erie~\eqref{eq:vasilyev} s'exprime aussi  
comme une s\'erie simple \`a laquelle le th\'eor\`eme~\ref{thwp} ci-dessus s'applique. 
Zlobin~\cite{Zlobincoeff} a r\'ecemment obtenu une d\'emonstration totalement diff\'erente 
par une \'etude directe de la s\'erie~\eqref{eq:vasilyev}, dans l'esprit des m\'ethodes combinatoires 
d\'evelopp\'ees dans cet article. On peut alors d\'emontrer des r\'esultats essentiellement de m\^eme
nature que ceux de~\cite{BR, RivoalCRAS}, ce qui renforce l'int\'er\^et pour des sommes multiples 
 sur des indices ordonn\'es.

\bigskip

Nous avons d\'emontr\'e dans \cite{CFRalgo} que toute s\'erie convergente de la forme 
\eqref{eq012}  s'\'ecrit comme combinaison 
lin\'eaire de polyz\^etas de poids au plus $pA$ et de profondeur au plus $p$ (et ce r\'esultat th\'eorique a \'et\'e obtenu, ind\'ependamment, par Zlobin \cite{ZlobinZametki2005}). 
En outre, nous avons pr\'esent\'e un algorithme, que nous avons impl\'ement\'e \cite{CFRweb}
en Pari, pour calculer explicitement une telle  combinaison lin\'eaire. Ceci 
nous a permis de d\'ecouvrir les propri\'et\'es de sym\'etrie que nous 
\'enon\c cons maintenant\footnote{Pour simplifier, nous ne d\'emontrons ici  
le th\'eor\`eme \ref{thprof2} que dans le cas o\`u $n$ est pair : voir la remarque~\ref{rqparite}. Cependant, il nous semble raisonnable d'esp\'erer que ce th\'eor\`eme soit vrai aussi quand $n$ est impair.}  
dans le cas particulier de la profondeur 2 : 

\bigskip

\begin{Th} \label{thprof2}
Soient  $n \geq 0$ et $A \geq 1$ des entiers, avec $n$ pair. 
Soit $P \in \Q[X_1,   X_2]$ un polyn\^ome en deux variables, de degr\'e $\leq A(n+1)-2$ par 
rapport \`a chacune d'elles, tel que
\begin{equation} \label{eqdefadeux}
\left\{
\begin{array}{l}
P(X_1, X_2) = - P(X_2, X_1) \\
P(-n-X_1, X_2 ) = (-1)^{A(n+1)+1} P( X_1, X_2 )\\
P(X_1, -n-X_2 ) = (-1)^{A(n+1)+1} P( X_1, X_2 )
\end{array}
\right.
\end{equation}
 Alors la somme double  \eqref{eq012} est une combinaison lin\'eaire, \`a coefficients rationnels :
 \begin{itemize}
 \item de 1, 
 \item de valeurs $\zeta(s)$ avec $s$ entier impair compris au sens large entre 3 et 2A,
 \item de diff\'erences $\zeta(s,s') - \zeta(s',s)$ avec $s$, $s'$ entiers impairs tels que $3 \leq s < s' \leq A$.
 \end{itemize} 
\end{Th}

Bien entendu, parmi les conditions \eqref{eqdefadeux}, la troisi\`eme est cons\'equence des deux premi\`eres. 
En particulier, si $A=4$, ce th\'eor\`eme montre que la s\'erie double
$$
\sum_{k_1 \geq   k_2 \geq 1} 
\frac{P(k_1,  k_2)}{(k_1)_{n+1}^{4}  (k_2)_{n+1}^{4} }
$$
est une forme lin\'eaire en $1$, $\zeta(3)$,  $\zeta(5)$ et  $\zeta(7)$ (ce qui \'etait 
loin d'\^etre \'evident a priori puisqu'on part d'une s\'erie double). Pour $A=3$, on obtient 
 une forme lin\'eaire en $1$, $\zeta(3)$,  $\zeta(5)$ ; enfin, pour $A=2$, une forme lin\'eaire en $1$ et $\zeta(3)$. 

\bigskip

Il est \`a noter que dans la s\'erie \eqref{eq012}, les variables $k_1$, \ldots, $k_p$ 
sont li\'ees par des in\'egalit\'es {\em larges}, comme dans \cite{CFRalgo} mais \`a l'inverse de la d\'efinition des polyz\^etas. 

\bigskip

Par exemple, le th\'eor\`eme \ref{thprof2} donne le cas particulier suivant : 

\begin{Cor} \label{corprofdeux}
Soient $n, r, t, \eps  \geq 0$ et $A \geq 1$ des entiers, avec $n$ pair,  tels que
$$\eps  \equiv (A+1)(n+1) +1  \bmod 2$$ 
et
$$\eps + 4 r + 2 t \leq (A-1)(n+1)-4.$$
Alors la s\'erie convergente
$$
\sum_{k_1 \geq k_2 \geq 1}
  \big(k_1 + \frac{n}{2}\big)^{\eps }  \big(k_2 + \frac{n}{2}\big)^{\eps } 
 \frac{ (k_1-k_2-r)_{2r+1} (k_1+k_2+n-r)_{2r+1} (k_1-t)_{2t+n+1}  
 (k_2-t)_{2t+n+1} }{(k_1)_{n+1}^{A}  \;(k_2)_{n+1}^{A} }
$$
est une combinaison lin\'eaire, \`a coefficients rationnels, de 1,   
de valeurs $\zeta(s)$ (avec $s$ entier impair tel que  $3 \leq s \leq  2A-1$), 
et  de diff\'erences $\zeta(s,s') - \zeta(s',s)$ (avec $s$, $s'$ entiers impairs tels que $3 \leq s < s' \leq A$).
\end{Cor}

Par exemple, on a 
\begin{multline*}
\sum_{k_1\ge k_2\ge 1} \big(k_1+\frac 12\big)\big(k_2+\frac 12\big)\frac{(k_1-k_2-1)_3(k_1+k_2)_3(k_1-1)_4(k_2-1)_4}
{(k_1)_2^7\;(k_2)_2^7} 
\\
= -1156 +891\,\zeta(3)+ \frac{189}2 \,\zeta(5) + 78 \big(\zeta(5,3) -\zeta(3,5)\big).
\end{multline*}

\bigskip

Un autre ingr\'edient, qui est fr\'equemment utilis\'e avec des s\'eries simples, consiste \`a 
d\'eriver la fraction rationnelle en $k$, avant de sommer ; par exemple, une double d\'erivation sert \`a montrer le r\'esultat de 
Zudilin~\cite{Zudilinonze} rappel\'e apr\`es la conjecture~\ref{conj:zeta alg inde}. 
Cette astuce, appliqu\'ee plusieurs fois, permet de faire dispara\^{\i}tre $\zeta(s)$ de la forme lin\'eaire
obtenue, pour de petites valeurs de $s$. On peut imaginer de l'utiliser pour des sommes multiples, m\^eme 
si on n'a aucun r\'esultat connu de disparition de polyz\^etas dans ce cadre (voir cependant \cite{SFHoffman}). Il est clair qu'en d\'erivant 
une fraction rationnelle de la forme 
$P(X_1, \ldots, X_p)/\big((X_1)_{n+1}^{A} \ldots (X_p)_{n+1}^{A}\big)$ par rapport \`a l'une des variables 
$X_i$, on obtient une fraction rationnelle de la m\^eme forme (avec $A$ remplac\'e par $A+1$). En 
profondeur 2, si un polyn\^ome $P(X_1, X_2)$ v\'erifie les relations \eqref{eqdefadeux}, alors 
le polyn\^ome $Q$ d\'efini par 
$$
\Big( \frac{\partial}{\partial X_1} \Big)^2  \Big(\frac{\partial}{\partial X_2} \Big)^2 
\frac{P(X_1, X_2)}{(X_1)_{n+1}^{A}  (X_2)_{n+1}^{A} }
 =   \frac{Q(X_1, X_2)}{(X_1)_{n+1}^{A+2}  (X_2)_{n+1}^{A+2} }
$$
 les v\'erifie aussi ; on peut donc lui appliquer aussi le th\'eor\`eme \ref{thprof2}. Cette 
remarque montre qu'on aurait pu ajouter des d\'erivations dans le corollaire \ref{corprofdeux}. 
Elle s'applique aussi en profondeur quelconque.
 
 \bigskip

Ce texte est divis\'e comme suit. 
Nous donnons au paragraphe \ref{sec2} l'\'enonc\'e g\'en\'eral, en profondeur quelconque, 
que nous obtenons. C'est l'occasion d'introduire la notion de {\em polyz\^etas antisym\'etriques}, et aussi de comparer notre g\'en\'eralisation des s\'eries (tr\`es) bien \'equilibr\'ees \`a celles provenant des syst\`emes de racines.

La preuve utilise deux outils : la r\'egularisation des s\'eries 
\`a divergence logarithmique et le d\'eveloppement en \'el\'ements simples des fractions 
rationnelles, qui sont pr\'esent\'es aux paragraphes \ref{sec3} et \ref{sec4} respectivement. 
Ces outils permettent d'\'enoncer (au paragraphe \ref{subsec42}) le th\'eor\`eme \ref{thconj3},
 qui implique notre r\'esultat principal (voir \S \ref{subsec43}). Ce th\'eor\`eme est 
d\'emontr\'e au paragraphe \ref{sec5}, par r\'ecurrence sur la profondeur : il s'agit 
du c{\oe}ur de la preuve. Le cas des profondeurs 1, 2 et 3 sont d\'etaill\'es s\'epar\'ement, et 
servent d'introduction \`a la d\'emonstration g\'en\'erale.

Enfin, au paragraphe \ref{sec6}, on d\'emontre le th\'eor\`eme \ref{thdecouple} \'enonc\'e ci-dessus. 
La preuve suit la m\^eme strat\'egie que celle du r\'esultat principal, mais chaque 
\'etape est nettement plus simple \`a mettre en {\oe}uvre. 

\bigskip

\noindent{\bf Remerciements : } Les auteurs ont eu l'opportunit\'e d'utiliser 
la puissance de calcul de la grappe M\'edicis, ce qui leur a permis de mener plus facilement les
exp\'erimentations qui ont conduit aux r\'esultats de cet article. Nous remercions \'egalement C.~Krattenthaler, M. Schlosser, W. Zudilin et l'arbitre 
pour leurs nombreuses remarques sur cet article, en particulier pour avoir port\'e \`a notre 
attention le lien entre nos s\'eries et les syst\`emes de racines. Enfin, le premier auteur remercie l'I.H.\'E.S. pour l'invitation lors de laquelle il a pu terminer ce travail.

\section{L'\'enonc\'e dans le cas convergent} \label{sec2}

\subsection{Polyz\^etas antisym\'etriques} \label{subsec21}

Pour \'enoncer notre r\'esultat en profondeur quelconque, nous aurons besoin de la 
notation suivante. Pour $p \geq 0$ et $s_1, \ldots, s_p \geq 2$ entiers, on pose
$$\zetaantisym (s_1, \ldots, s_p) =  \sum_{\sigma\in\mathfrak{S}_p} \eps_{\sigma} 
\zeta(s_{\sigma(1)}, \ldots, s_{\sigma(p)}) ,$$
o\`u $\eps_{\sigma} $ d\'esigne la signature de la permutation $\sigma$. 
On appelle {\em polyz\^eta antisym\'etrique} une telle combinaison lin\'eaire de polyz\^etas 
(m\^eme si, pour $p \geq 2$, ce n'est pas en g\'en\'eral un polyz\^eta). Il s'agit de s\'eries 
convergentes, puisque tous les $s_i$ sont suppos\'es \^etre sup\'erieurs ou \'egaux \`a 2 ; on 
utilisera donc parfois le terme de {\em polyz\^eta antisym\'etrique convergent}. Pour $p=1$, on 
a $\zetaantisym (s) = \zeta(s)$. La convention naturelle consiste \`a poser $\zetaantisym (s_1, 
\ldots, s_p) = 1$ lorsque $p=0$, puisqu'il existe une unique bijection de l'ensemble vide dans lui-m\^eme.  Pour $p=2$, on a 
$\zetaantisym (s_1, s_2) = \zeta(s_1, s_2) - \zeta(s_2, s_1)$ 
et lorsque $p=3$, on a 
\begin{multline*}
\zetaantisym (s_1, s_2, s_3) 
\\
= \zeta(s_1, s_2, s_3) + \zeta(s_2, s_3, s_1) +  \zeta(s_3, s_1, s_2)
- \zeta(s_2, s_1, s_3) - \zeta(s_1, s_3, s_2) -  \zeta(s_3, s_2, s_1).
\end{multline*} 
Par d\'efinition, pour tout $\sigma\in\mathfrak{S}_p$ on a 
$$ \zetaantisym(s_{\sigma(1)}, \ldots, s_{\sigma(p)}) =  \eps_{\sigma} \zetaantisym (s_1, \ldots, s_p), $$
et $\zetaantisym (s_1, \ldots, s_p) = 0$ d\`es que deux des $s_i$ sont \'egaux. 

\bigskip

Il nous semble raisonnable de penser qu'en g\'en\'eral, un poyz\^eta antisym\'etrique n'est pas un polyn\^ome en valeurs de la fonction $\zeta$ de Riemann. En revanche, tout polyz\^eta ``sym\'etrique'' (d\'efini comme $ \zetaantisym (s_1, \ldots, s_p)$ mais en omettant la signature $ \eps_{\sigma} $)
est un polyn\^ome en les valeurs $\zeta(s)$ (d'apr\`es \cite{Hoffman1992}, Theorem 2.2). 

\subsection{Enonc\'e du r\'esultat principal} \label{subsec22}

Notons   $\mathscr{A}_p$   l'ensemble des polyn\^omes 
$P(X_1, \ldots, X_p)\in\Q[X_1, \ldots, X_p]$  tels que :
$$
\begin{cases}
\mbox{Pour tout } \sigma \in\mathfrak{S}_p \mbox{,  on ait }  
\\
\qquad \qquad 
P(X_{\sigma(1)}, X_{\sigma(2)},\ldots, X_{\sigma(p)}) = \eps_{\sigma} P(X_1, X_2, \ldots, X_p). 
\\
\\
\mbox{Pour tout } j \in \unp \mbox{,  on ait}
\\ 
\qquad \qquad  P(X_1,\ldots, X_{j-1}, -X_j-n, X_{j+1}, \ldots, X_p ) 
\\
\qquad \qquad \qquad \qquad = (-1)^{A(n+1)+1} P( X_1,\ldots, X_{j-1}, X_j, X_{j+1}, \ldots, X_p ).
\end{cases}
$$
Ces conditions (qui font appara\^{\i}tre l'action de groupe qui sera utilis\'ee au paragraphe \ref{subsec41}) sont bien s\^ur redondantes. Si la premi\`ere est satisfaite, alors il suffit notamment  de v\'erifier la seconde pour une seule valeur  de $j$.

Par exemple, $\mathscr{A}_2$ est exactement l'ensemble des polyn\^omes $P$  v\'erifiant les conditions \eqref{eqdefadeux}. 
Par ailleurs, si $P \in \mathscr{A}_p$ alors $P$ a le m\^eme degr\'e par rapport \`a chacune des variables 
$X_1,   \ldots, X_p$. Bien entendu la d\'efinition de $\mathscr{A}_p$ d\'epend aussi de la 
parit\'e de $A(n+1)$, mais on ne refl\`ete pas cette d\'ependance pour ne pas alourdir la notation. 

\bigskip

Nous pouvons maintenant \'enoncer notre r\'esultat principal.\footnote{Ce r\'esultat, comme 
les th\'eor\`emes \ref{thconj2} et \ref{thconj3} ci-dessous, ne sera d\'emontr\'e ici que dans 
le cas o\`u  $n$ est pair. Ceci permet de simplifier la preuve (voir la remar\-que~\ref{rqparite})  
et ne devrait pas \^etre un obstacle \`a d'\'eventuelles applications diophantiennes. Cependant,  il nous semble raisonnable d'esp\'erer que ces \'enonc\'es soient vrais aussi quand $n$ est impair.}

\begin{Th} \label{thconj1}
Soient $n \geq 0$ et $A, p \geq 1$ des entiers, avec $n$ pair. 
Soit $P \in \mathscr{A}_p$ de degr\'e $\leq A(n+1)-2$ par rapport \`a chacune des variables. Alors la s\'erie
\begin{equation} \label{eq021}
\sum_{k_1 \geq  \ldots \geq   k_p \ge 1} 
\frac{P(k_1, \ldots, k_p)}{(k_1)_{n+1}^{A} \ldots (k_p)_{n+1}^{A} }
\end{equation}
est une combinaison lin\'eaire, \`a coefficients rationnels, de produits de la forme
$$\zeta(s_1) \ldots  \zeta(s_{q}) \zetaantisym (s'_{1}, \ldots, s'_{q'})$$
avec 
\begin{equation} \label{eqcondithconj1}
\left\{
\begin{array}{l}
q, q' \geq 0   \mbox{ entiers tels que } 2q+q' \leq p, \\
s_1, \ldots, s_q, s'_1, \ldots, s'_{q'} \mbox{ entiers impairs } \geq 3, \\
s_i \leq 2A-1  \mbox{ pour tout } i \in \unq, \\
s'_i \leq A  \mbox{ pour tout } i \in \unqpr.
\end{array}
\right.
\end{equation}
\end{Th}

La dissym\'etrie entre $s_1, \ldots, s_q$ d'une part, et $ s'_1, \ldots, s'_{q'} $ 
d'autre part, dans la conclusion de cet \'enonc\'e sera comment\'ee plus loin 
(juste apr\`es l'\'enonc\'e du th\'eor\`eme \ref{thconj3}).

\smallskip

Il est important de bien visualiser l'ensemble des produits de polyz\^etas qui apparaissent dans ce th\'eor\`eme. Par exemple, lorsque $q'=0$ le polyz\^eta antisym\'etrique $ \zetaantisym (s'_{1}, \ldots, s'_{q'})$ vaut 1 (conform\'ement \`a la convention \'evoqu\'ee au paragraphe \ref{subsec21}), et on obtient un produit de valeurs de $\zeta$ en des entiers impairs. Lorsque $q=q'=0$, ce produit est vide et on obtient~1.

Si $p=1$, le th\'eor\`eme \ref{thconj1} affirme que \eqref{eq021} est une
 combinaison lin\'eaire de $1$ et des $\zeta(s)$ pour $s$ impair tel que $3 \leq s \leq A$ : 
on retrouve le th\'eor\`eme \ref{thwp}, c'est-\`a-dire le ph\'enom\`ene de sym\'etrie  
li\'e aux s\'eries hyperg\'eom\'etriques (tr\`es) bien \'equilibr\'ees  en profondeur $1$.

Si $p=2$, on obtient exactement le th\'eor\`eme \ref{thprof2} \'enonc\'e dans l'introduction.

Si $p=3$, ce th\'eor\`eme affirme que la s\'erie est une combinaison lin\'eaire, \`a coefficients rationnels :
 \begin{itemize}
 \item de produits d'au plus deux valeurs  de $\zeta$ en des entiers impairs $\geq 3$, 
 \item de polyz\^etas antisym\'etriques convergents $\zetaantisym (s_1, s_2)$ avec $s_1, s_2 \geq 3$ impairs, 
 \item de polyz\^etas antisym\'etriques convergents $\zetaantisym (s_1, s_2, s_3)$ avec $s_1, s_2, s_3 \geq 3$ impairs. 
 \end{itemize} 

En profondeur $p \geq 4$, des termes tels que $q \geq 1$ et $q' \geq 2$ peuvent appara\^{\i}tre : 
il semble que la s\'erie obtenue ne soit pas toujours la somme d'un polyn\^ome en valeurs  $\zeta(s)$ 
(avec $s $ impair) et d'une combinaison lin\'eaire de polyz\^etas antisym\'etriques $\zetaantisym 
(s_1, \ldots, s_q)$ avec $s_1, \ldots, s_q $ impairs. 

\`A l'inverse, on peut affaiblir la conclusion du th\'eor\`eme \ref{thconj1} en disant que la
 s\'erie est un polyn\^ome (\`a coefficients rationnels) en les polyz\^etas antisym\'etriques 
convergents $\zetaantisym (s_1, \ldots, s_q)$ avec $1 \leq q \leq p$ et $s_1, \ldots , s_q \geq 3$ 
impairs tels que $s_1 +\ldots + s_q \leq pA$.

\bigskip

Lorsque $A \leq 2$,  on a forc\'ement $q'=0$ pour tous les produits qui apparaissent, ce qui fournit le corollaire suivant :

\begin{Cor} 
Sous les hypoth\`eses du th\'eor\`eme \ref{thconj1}, si $A \leq 2$ alors la s\'erie \eqref{eq021} est 
 un polyn\^ome en $\zeta(3)$ \`a coefficients rationnels. 
\end{Cor}

\bigskip

Le th\'eor\`eme \ref{thconj1} contient, par exemple, le cas particulier suivant :

\begin{Cor} \label{corconj1}
Soient $n, r, t, \eps \geq 0$ et $A,p \geq 1$ des entiers, avec $n$ pair,  tels que
$$\eps  \equiv (A+1)(n+1) +1  \bmod 2$$
et
$$\eps + (4r+2)p + 2t \leq (A-1)(n+1) + 4r.$$
Alors la s\'erie convergente
\begin{equation}\label{eq:seriecorollaire}
\sum_{k_1 \geq \ldots \geq  k_p \geq 1}
\bigg[\prod_{i=1} ^p (k_i + \frac{n}{2})\bigg]^{\eps} 
 \frac{ \displaystyle  \bigg[ \prod_{1 \leq i < j \leq p}  (k_i-k_j-r)_{2r+1} (k_i+k_j+n-r)_{2r+1} \bigg] 
 \bigg[ \prod_{i=1} ^p (k_i-t)_{2t+n+1}   \bigg]  }{(k_1)_{n+1}^{A} \ldots  (k_p)_{n+1}^{A} }
\end{equation}
est une combinaison lin\'eaire comme celles du th\'eor\`eme \ref{thconj1}. 
\end{Cor}

Un  exemple d'application de ce corollaire est la s\'erie suivante (dans laquelle on prend $t=0$ et les symboles de Pochhammer $(k_i)_{n+1}$ se simplifient avec ceux du d\'enominateur) : 
\begin{multline}\label{eq:exempleinteressant}
\sum_{k_1\ge k_2\ge k_3\ge 1}\big(k_1+\frac12\big)\big(k_2+ \frac12\big)\big(k_3+\frac12\big) 
\\
\times \frac{(k_1-k_2)(k_2-k_3)(k_1-k_3)
(k_1+k_2+1)(k_1+k_3+1)(k_2+k_3+1)}{(k_1)_2^4\;(k_2)_2^4\;(k_3)_2^4}
\\
= -\frac{1}{4} - \zeta(3) + \frac14 \,\zeta(5) + \zeta(3)^2 -\frac14 \,\zeta(7).
\end{multline}

\medskip

Dans d'\'eventuelles applications diophantiennes (voir \S \ref{subsec23}), on pourrait prendre $\eps$   \'egal \`a 0 ou 1, de telle sorte que sa contribution asymptotique (pour $n$ grand) serait n\'egligeable. Le probl\`eme est de bien choisir les param\`etres $r$ et $s$ en fonction de $n$, ou encore d'imaginer d'autres 
polyn\^omes $P$ auxquels on pourrait appliquer le th\'eor\`eme \ref{thconj1}.

\bigskip

On pourrait chercher \`a obtenir un analogue du th\'eor\`eme \ref{thconj1} dans lequel seuls des entiers $s_i$ et $s'_i$ {\em pairs} appara\^{\i}traient. Un tel \'enonc\'e correspondrait peut-\^etre \`a des polyn\^omes $P$ invariants sous l'action de $\spp$, \`a des polyz\^etas {\em sym\'etriques} (voir la fin du paragraphe \ref{subsec21}), ou \`a des valeurs de polylogarithmes en un point $z  = -1$ (c'est-\`a-dire \`a un signe, d\'ependant de $k_1$, \ldots, $k_p$,  qui multiplierait la fraction rationnelle que l'on somme). 

\bigskip

Toujours en vue d'une \'eventuelle application diophantienne, 
il serait utile d'avoir 
un contr\^ole sur le d\'enominateur des coefficients qui interviennent dans l'\'ecriture de \eqref{eq021} comme combinaison lin\'eaire de polyz\^etas. Lorsque  $P = n!^{Ap} \widetilde P $ o\`u $\widetilde P $ est un polyn\^ome \`a coefficients entiers, on peut supposer dans le th\'eor\`eme \ref{thconj1} que  $\dd_n ^{Ap}$ est un d\'enominateur commun des coefficients de la combinaison lin\'eaire (o\`u $\dd_n$ est le ppcm des entiers 1, 2, \ldots, $n$ ; ceci sera d\'emontr\'e au paragraphe \ref{subsecdenompreuve}). Dans certains autres cas, la pr\'esence de symboles de Pochhammer dans la d\'efinition de $P$ permet d'obtenir un tel d\'enominateur, comme c'est le cas habituellement en profondeur 1. \'Etant donn\'e un polyn\^ome $P$ particulier, il n'est pas difficile  de d\'eduire un tel r\'esultat  du th\'eor\`eme  \ref{thconj3} ci-dessous (il suffit d'adapter le lemme \ref{lemdenom} qui figure au paragraphe \ref{subsecdenompreuve}). En outre, il serait int\'eressant de savoir si une {\em conjecture des d\'enominateurs} analogue \`a celle d\'emontr\'ee dans \cite{KR} existe.
 
\subsection{Liens avec les s\'eries hyperg\'eom\'etriques issues de syst\`emes de racines}
			 
Lorsque l'on ne pr\'ecise pas la forme du polyn\^ome 
$P(X_1,X_2,\ldots, X_p)\in\mathbb{Q}[X_1, X_2, \ldots, X_p]$ au 
num\'erateur de \eqref{eq012}, 
nos s\'eries multiples peuvent s'exprimer comme combinaisons
lin\'eaires \`a coefficients rationnels de s\'eries hyperg\'eom\'etriques multiples 
de Lauricella. Lorsque $p=1$, la s\'erie~\eqref{eq:seriecorollaire} consid\'er\'ee au corollaire~\ref{corconj1} est une s\'erie simple hyperg\'eom\'etrique {\em very well-poised}. 

Il est donc naturel de se demander si, pour $p \ge 2$, la s\'erie multiple~\eqref{eq:seriecorollaire} correspond \`a l'une ou l'autre des g\'en\'eralisations de {\em well-poisedness} en dimension sup\'erieure, qui sont li\'ees aux syst\`emes de racines $C_n$, $D_n$ ou $BC_n$ 
(voir par exemple~\cite{humphreys} pour les d\'efinitions). On peut faire les remarques suivantes.
Dans~\cite{schlosser}, une s\'erie hyperg\'eom\'etrique multiple est dite 
de type $C_n$ si le facteur 
\begin{equation}\label{eq:cn}
\Big(\prod_{1\le i<j\le n} (k_i-k_j+x_i-x_j)(k_i+k_j+x_i+x_j) \Big)
\Big(\prod_{i=1}^n (k_i+x_i)\Big)
\end{equation}
est pr\'esent, la sommation \'etant sur les $k_1 \ge 0, k_2\ge 0, \ldots, k_n\ge 0$, les $x_j$ \'etant des param\`etres. Elle est dite de type $D_n$ si le facteur 
\begin{equation}\label{eq:dn}
\prod_{1\le i<j\le n} (k_i-k_j+x_i-x_j)(k_i+k_j+x_i+x_j)
\end{equation}
est pr\'esent mais pas le facteur $\prod_{i=1}^n (k_i+x_i)$. Le type $C_n$ est donc une des g\'en\'eralisations possibles des s\'eries  {\em very well-poised}, tandis que le type $D_n$ g\'en\'eralise les s\'eries qui sont   {\em well-poised} mais pas  {\em very well-poised}. Cependant, aucune de ces d\'efinitions n'impose de propri\'et\'e de sym\'etrie sur le sommande, alors que dans tous les \'enonc\'es obtenus ici les propri\'et\'es de sym\'etrie sont cruciales : des exemples (faciles \`a calculer gr\^ace \`a \cite{CFRweb}) permettent facilement de voir qu'on ne peut pas remplacer, dans nos r\'esultats, l'hypoth\`ese de sym\'etrie par une hypoth\`ese de divisibilit\'e par un facteur du type \eqref{eq:cn} ou \eqref{eq:dn}. 

Par exemple, dans le corollaire~\ref{corconj1} ci-dessus, pour $r=0$, le terme de
la s\'erie $p$-uple~\eqref{eq:seriecorollaire} est de type $C_p$ lorsque $\varepsilon=1$ et de type $D_p$ lorsque $\varepsilon=0$, avec $x_i=n/2+1$. La s\'erie triple~\eqref{eq:exempleinteressant} est, quant 
\`a elle, de type $C_p$, avec $x_i=3/2$.
Cependant, dans ces deux cas, notre sommation porte sur 
$k_1 \ge k_2\ge \cdots \ge k_p\ge 1$ ce qui, comme on va maintenant le voir, produit 
une tr\`es grosse diff\'erence sur la nature des polyz\^etas qui apparaissent. En effet, en modifiant la sommation dans~\eqref{eq:exempleinteressant}, on obtient l'\'evaluation d'une s\'erie de type~$C_3$~:
\begin{multline*}\label{eq:exempleinteressantbis}
\sum_{k_1, k_2, k_3\ge 1}\big(k_1+\frac12\big)\big(k_2+ \frac12\big)\big(k_3+\frac12\big) 
\\
\times \frac{(k_1-k_2)(k_2-k_3)(k_1-k_3)
(k_1+k_2+1)(k_1+k_3+1)(k_2+k_3+1)}{(k_1)_2^4\;(k_2)_2^4\;(k_3)_2^4}
= 0
\end{multline*}
puisque le sommande est chang\'e en son oppos\'e par l'\'echange des indices $k_1 \leftrightarrow k_2$. Cette remarque vaut aussi pour la somme de type $D_3$ :
\begin{equation*}
\sum_{k_1, k_2, k_3\ge 1}
\frac{(k_1-k_2)(k_2-k_3)(k_1-k_3)
(k_1+k_2+1)(k_1+k_3+1)(k_2+k_3+1)}{(k_1)_2^4\;(k_2)_2^4\;(k_3)_2^4}
= 0.
\end{equation*}

Le choix de l'ensemble de sommation des s\'eries est donc crucial afin d'obtenir des r\'esultats non 
triviaux \`a partir de  s\'eries pr\'esentant les sym\'etries $C_n$ et $D_n$. Par ailleurs, on peut remarquer que ces deux sym\'etries ne tiennent finalement que tr\`es peu compte de la forme des 
sommandes des s\'eries telles que~\eqref{eq:seriecorollaire}. 
Michael Schlosser nous a fait remarquer que ces s\'eries pr\'esentent en fait 
une sym\'etrie proche du type $BC_n$, qui tient compte de la pr\'esence de facteurs 
{\og Pochhammer\fg} et dont l'\'etude est toute r\'ecente (voir~\cite{coskun}). Les sym\'etries issues des divers syst\`emes de racines ont donc un grand int\'er\^et dans l'\'etude diophantienne des  polyz\^etas et on peut esp\'erer qu'elles puissent jouer un r\^ole de plus en plus  important \`a l'avenir.

\subsection{Applications diophantiennes \'eventuelles} \label{subsec23}
 
Pour tout entier $A \geq 1$, notons $\calf_A$ le sous-$\Q$-espace vectoriel de $\R$ engendr\'e par 1 et les $\zeta(s)$, pour $s$ entier impair tel que $3 \leq s \leq A$. Les minorations suivantes sont essentiellement les seules connues (voir par exemple \cite{SFBou} pour un survol) : 
\begin{equation} \label{eqminodim}
\begin{cases}
\dim \calf_3 = 2 \quad \mbox{ \cite{Apery}}
\\
\dim \calf_{145} \geq 3 \quad \mbox{ (\cite{Zudilincentqc}, voir aussi \cite{BR})} 
\\
\dim \calf_A \geq \frac{1 - \odu}{1 + \log 2} \log A \quad \mbox{ (\cite{BR}, \cite{RivoalCRAS})}.
\end{cases}
\end{equation}

\bigskip

Pour $A \geq 1$ et $p \geq 1$, notons $\cale_{A, p}$ le sous-$\Q$-espace vectoriel de $\R$ engendr\'e par les produits $\zeta(s_1) \ldots  \zeta(s_{q}) \zetaantisym (s'_{1}, \ldots, s'_{q'})$ satisfaisant aux conditions \eqref{eqcondithconj1} \'enonc\'ees dans le th\'eor\`eme \ref{thconj1}. L'int\'er\^et de ce th\'eor\`eme est justement de fournir des s\'eries qui appartiennent \`a $\cale_{A, p}$, et qui pourraient permettre de minorer la dimension de cet espace.

Pour $p= 1$ on a simplement $\cale_{A, 1} = \calf_A$. Pour $p \geq 2$, l'inclusion $\calf_{2A-1} \subset \cale_{A, p}$ permet d'obtenir, \`a partir de \eqref{eqminodim}, des minorations de $\dim \cale_{A, p}$. On peut esp\'erer que le th\'eor\`eme \ref{thconj1} (ou le corollaire \ref{corconj1}) conduisent \`a des minorations plus fines de  $\dim \cale_{A, p}$, qui constitueraient de nouveaux r\'esultats diophantiens. Par exemple,  peut-\^etre peut-on obtenir une minoration de la forme $\dim \cale_{A, p} \geq (c(p)-\opdu) \log(A)$, o\`u $\opdu$ est une suite qui d\'epend de $p$ et $A$ et tend vers 0 quand $A$ tend vers l'infini (quelle que soit la valeur, fix\'ee, de $p$), et $c(p)$ est une fonction de $p$ seulement. Ceci serait nouveau \`a condition qu'on ait $c(p) > \frac{1}{1 + \log 2}$ (ce que l'on peut esp\'erer, notamment si $p$ est grand). 

\bigskip

Par ailleurs, si on arrivait \`a montrer que $\dim \cale_{2, p} \geq 3$ pour un certain $p$, on obtiendrait que $\zeta(3)$  n'est pas quadratique. Si cette dimension pouvait \^etre arbitrairement grande, cela donnerait la transcendance de $\zeta(3)$. Malheureusement, les contraintes de sym\'etrie impos\'ees au polyn\^ome $P$ dans le th\'eor\`eme \ref{thconj1} semblent trop draconiennes pour qu'on puisse aboutir \`a un r\'esultat aussi spectaculaire (voir \`a ce propos \cite{SFHoffman}, o\`u des propri\'et\'es de sym\'etrie plus faibles sont d\'emontr\'ees sous des hypoth\`eses moins restrictives). Cependant, l'une des motivations principales de cet article est de montrer que l'algorithme de \cite{CFRalgo} permet de deviner des propri\'et\'es, comme celles d\'emontr\'ees ici, de disparition de polyz\^etas. La structure de la preuve du th\'eor\`eme \ref{thconj1} devrait pouvoir \^etre utilis\'ee pour d\'emontrer d'autres r\'esultats analogues, dont les applications diophantiennes pourraient \^etre plus faciles.

\section{R\'egularisation des s\'eries divergentes} \label{sec3}

\subsection{Rappels} \label{subsec31}

Dans toute la suite, on note $H_N$ la somme harmonique d\'efinie par 
$$H_N = 1 + \frac12 + \frac13 + \ldots + \frac{1}{N}.$$
La proposition suivante a \'et\'e d\'emontr\'ee par   Racinet (voir le Corollaire 2.1.8 
de \cite{RacinetIHES}), en suivant des travaux de Boutet de Monvel.

\begin{Prop} \label{propstuffle} 
Soient $p \geq 0$ et $s_1, \ldots, s_p \geq 1$. Alors il existe un unique polyn\^ome $Q$ tel que, 
  pour tout $\eps > 0$, on ait quand $N$ tend vers $+\infty$ : 
$$ \sum_{N \geq k_1 > \ldots > k_p \geq 1} 
\frac{1}{k_1 ^{s_1} \ldots k_p ^{s_p}} = Q(H_N) + \gdo_\eps (N^{-1+\eps});
$$
on note alors $\zetaetoile (s_1,\ldots, s_p)$ le coefficient constant de $Q$, c'est-\`a-dire sa valeur en 0.
\end{Prop}

Cette proposition d\'efinit les {\em valeurs r\'egularis\'ees} $\zetaetoile (s_1,\ldots, s_p)$ des s\'eries divergentes 
$$\sum_{ k_1 > \ldots >  k_p \geq 1} 
\frac{1}{k_1 ^{s_1} \ldots k_p ^{s_p}}$$
lorsque $s_1 =1$. D\`es que $s_1 \geq 2$, on a simplement $\zetaetoile (s_1,\ldots, s_p) = \zeta (s_1,\ldots, s_p)$ et le polyn\^ome $Q$ est constant. 

Il s'agit de la r\'egularisation relative au produit nomm\'e {\em stuffle} (voir \cite{MiW}), avec la convention $\zetaetoile (1) = 0$. Il existe une autre forme de r\'egularisation, li\'ee au produit shuffle, et utilis\'ee dans \cite{CFRalgo} ; mais nous n'en aurons pas besoin ici. 

\bigskip

Les valeurs r\'egularis\'ees $\zetaetoile (s_1,\ldots, s_p)$  peuvent se calculer de mani\`ere algorithmique ; ce sont des combinaisons lin\'eaires \`a coefficients rationnels de polyz\^etas. 

\bigskip

Nous aurons aussi besoin de la d\'efinition suivante. On appelle {\em polyz\^eta antisym\'etrique r\'egularis\'e} la combinaison  lin\'eaire suivante de polyz\^etas  r\'egularis\'es, pour  $p \geq 1$ et $s_1, \ldots, s_p \geq 1$ entiers :
$$\zetaetoileantisym (s_1, \ldots, s_p) =  \sum_{\sigma\in\mathfrak{S}_p} \eps_{\sigma} \zetaetoile(s_{\sigma(1)}, \ldots, s_{\sigma(p)}) .$$
Lorsque $s_1 \geq 2$, on a $\zetaetoileantisym (s_1, \ldots, s_p) = \zetaantisym (s_1, \ldots, s_p) $ : on retrouve les polyz\^etas antisym\'etriques convergents.  Lorsque $p=0$, on pose $\zetaetoileantisym (s_1, \ldots, s_p) = \zeta  (s_1, \ldots, s_p) = 1$.

\subsection{\'Enonc\'e avec r\'egularisation des divergences} \label{subsec32}

L'une des motivations principales pour consid\'erer des polyz\^etas r\'egularis\'es est qu'ils permettent de rendre la th\'eorie  plus compl\`ete, et en tout cas plus \'el\'egante. Nous en donnons ici une illustration : pour d\'emontrer le   th\'eor\`eme \ref{thconj1} (qui concerne seulement des s\'eries convergentes), nous allons utiliser le r\'esultat suivant (dans lequel des divergences logarithmiques sont autoris\'ees, et r\'egularis\'ees).\footnote{Plus pr\'ecis\'ement, nous d\'emontrerons au \S \ref{sec5} le th\'eor\`eme \ref{thconj3}, qui est une forme plus pr\'ecise  du  th\'eor\`eme \ref{thconj2}, et  nous en d\'eduirons le  th\'eor\`eme \ref{thconj1}   au paragraphe \ref{subsec43}.}

\begin{Th} \label{thconj2}
Supposons $n$ pair. 
Soit $P \in \mathscr{A}_p$ de degr\'e $\leq A(n+1)-1$ par rapport \`a chacune des variables. Alors il existe un polyn\^ome $Q_P$ tel que, pour tout $\eps > 0$, on ait quand $N$ tend vers $+\infty$ : 
\begin{equation} \label{eq1}
 \sum_{N \geq k_1 \geq \ldots \geq k_p \geq 1} 
\frac{P(k_1,\ldots,k_p)}{(k_1)_{n+1}^A \ldots (k_p)_{n+1}^A}
= Q_P (H_N) + \gdo_\eps (N^{-1+\eps}),
\end{equation}
et tel que $Q_P(0)$ soit  une combinaison lin\'eaire, \`a coefficients rationnels,
de produits de la forme
\begin{equation} \label{eq998}
\zetaetoile(s_1) \ldots  \zetaetoile(s_{q}) \zetaetoileantisym (s'_{1}, \ldots, s'_{q'})
\end{equation}
avec 
$$
\left\{
\begin{array}{l}
q, q' \geq 0 \mbox{ entiers tels que } 2q+q' \leq p \\
s_1, \ldots, s_q , s'_1, \ldots, s'_{q'} \mbox{ entiers impairs } \geq 1 \\
s_i \leq 2A-1  \mbox{ pour tout } i \in \unq \\
s'_i \leq A  \mbox{ pour tout }  i \in \unqpr .
\end{array}
\right.
$$
\end{Th}

Comme $\zetaetoile(1) = 0$, on peut se restreindre aux produits \eqref{eq998} tels que $s_1, \ldots, s_{q} \geq 3$. 

\bigskip

Si dans ce th\'eor\`eme on suppose que $P$ est de degr\'e $\leq A(n+1)-2$ par rapport \`a chacune des variables, alors \eqref{eq1} converge quand $N$ tend vers $+\infty$, donc le polyn\^ome $Q_P$ est constant (\'egal \`a $Q_P(0)$). Pour d\'eduire le th\'eor\`eme \ref{thconj1} du th\'eor\`eme \ref{thconj2}, il suffit donc de d\'emontrer que le produit \eqref{eq998} ne peut appara\^{\i}tre que si $ s'_1, \ldots, s'_{q'} \geq 3$.  C'est l'objet du paragraphe \ref{subsec43}; pour y parvenir, on utilise en fait une version plus pr\'ecise  du  th\'eor\`eme \ref{thconj2}, que nous allons formuler gr\^ace au d\'eveloppement en \'el\'ements simples.

\section{D\'ecomposition en \'el\'ements simples} \label{sec4}

\subsection{Notations et actions de groupes}  \label{subsec41}

Soit $P(k_1,\ldots,k_p)$ un polyn\^ome de degr\'e $\leq A(n+1)-1$ par rapport \`a chacune des variables, \`a coefficients rationnels. La d\'ecomposition en \'el\'ements simples de la fraction rationnelle 
\begin{equation} \label{eq041}
R(k_1,\ldots, k_p) = \frac{P(k_1,\ldots, k_p)}{(k_1)_{n+1}^A \ldots (k_p)_{n+1}^A}
\end{equation}
s'\'ecrit
\begin{equation} \label{eq2}
R(k_1,\ldots, k_p) = \sum_\indso{0 \leq j_1, \ldots, j_p \leq n}{1 \leq s_1, \ldots, s_p \leq A} 
\frac{\cjs}{(k_1+j_1)^{s_1} \ldots (k_p+j_p)^{s_p}}
\end{equation}
avec des rationnels $\cjs$. L'unicit\'e de ce d\'eveloppement montre que $P$ appartient \`a $\mathscr{A}_p$ si, et seulement si,  on a : 
\begin{equation} \label{eq3}
\left\{ 
\begin{array}{l}
\cjs = (-1)^{s_i+1} \cjsi \mbox{ pour tout } i \in \unp  \\
\\
\cjs = \eps_\gamma \cjsgamma \mbox{ pour tout } \gamma \in \spp.
\end{array}
\right.
\end{equation}

\bigskip

Donnons maintenant une interpr\'etation alg\'ebrique (en termes de groupes op\'erant sur des 
ensembles) de cette situation, qui sera utile dans les preuves.

\bigskip

Pour $\eps \in \zdz$ (o\`u on voit toujours $\zdz$ comme \'etant le groupe multiplicatif $\{-1, 1\}$) et $j \in \zeron$, on pose :
$$
\left\{
\begin{array}{l}
\eps \cdot j =j \mbox{ si } \eps =1 , \\
\eps \cdot j =n-j \mbox{ si } \eps = -1. 
\end{array}
\right.
$$
Ceci d\'efinit une action de $\zdz$ sur $\zeron$. De mani\`ere diagonale, on peut alors d\'efinir une action de $\zdzp$ sur $\zeron^p$ en posant :
$$(\eps_1,\ldots, \eps_p) \cdot (j_1,\ldots, j_p) = (\eps_1  \cdot j_1, \ldots,  \eps_p \cdot j_p).$$
En outre,  on consid\`ere l'action triviale de $\zdzp$ sur $\una^p$, et on en d\'eduit une action de $\zdzp$ sur $\zeron^p \times \una^p$ d\'efinie par :
$$(\eps_1,\ldots, \eps_p) \cdot (j_1,\ldots, j_p, s_1,\ldots, s_p) = (\eps_1  \cdot j_1, \ldots,  \eps_p \cdot j_p, s_1,\ldots, s_p).$$
Par ailleurs, le groupe $\spp$ agit par permutation des facteurs sur $\zdzp$,  sur  $\zeron^p$ et sur $\una^p$ (donc agit aussi sur $\zeron^p \times \una^p$). 
On en d\'eduit une action du produit semi-direct $\zdzp \psd \spp$ sur $\zeron^p \times \una^p$; et 
\eqref{eq3} signifie que $\cjs$ est constant (au signe pr\`es) sur chaque orbite (et ce signe est bien d\'etermin\'e en fonction de la parit\'e des $s_i$). 

\begin{remark} Le sous-groupe d'indice 2 de $\zdzp \psd \spp$ form\'e par les $(\eps_1,\ldots, \eps_p, 
\gamma)$ tel que $\eps_1 \ldots \eps_p = 1$ est d'ordre $2^{p-1} p!$ ; pour $p=5$, c'est exactement le 
groupe de Rhin-Viola \cite{RV3} pour $\zeta(3)$.  
Nous n'avons trouv\'e aucune explication \`a cette co\"{\i}ncidence.
\end{remark}

\subsection{\'Enonc\'e r\'egularis\'e en termes d'\'el\'ements simples}  \label{subsec42}

On va d\'eduire les th\'eor\`emes \ref{thconj1} et  \ref{thconj2} du r\'esultat suivant : 

\begin{Th} \label{thconj3}
Supposons $n$ pair. 
 Soient $j_1, \ldots, j_p \in \zeron$ et $s_1,\ldots,s_p  \geq 1$. Alors il existe un polyn\^ome $\qjs$ 
tel que,   pour tout $\eps > 0$, on ait quand $N$ tend vers $+\infty$ : 
\begin{multline} 
\sum_{N \geq k_1 \geq \ldots \geq k_p \geq 1} 
 \sum_{\sigma \in \spp}
 \sum_{(\eps_1,\ldots,\eps_p)\in \zdzp} 
\eps_\sigma  \eps_1^{s_1+1} \ldots \eps_p ^{s_p+1} 
\frac{1}{(k_{\sigma(1)} + \eps_1 \cdot j_1)^{s_1} \ldots (k_{\sigma(p)} + \eps_p \cdot j_p)^{s_p}} 
\\
= \qjs (H_N) + \gdo_\eps (N^{-1+\eps}), 
\label{eq4} 
\end{multline}
et tel que $\qjs(0)$ soit  une combinaison lin\'eaire, \`a coefficients rationnels, 
de produits de la forme
\begin{equation} \label{eq999}
\zetaetoile(s'_1) \ldots  \zetaetoile(s'_{q'}) \zetaetoileantisym (s''_{1}, \ldots, s''_{q''})
\end{equation}
avec, pour chaque produit de cette forme :  
$$
\left\{
\begin{array}{l}
q' ,q'' \geq 0   \mbox{ entiers tels que } 2q'+q'' \leq p \\
s'_1, \ldots, s'_{q'}, s''_{1}, \ldots, s''_{q''} \geq 1 \mbox{  impairs } \\
\mbox{il existe } \sigma \in \spp \mbox{ tel que : }\\
\quad \quad \bullet \, \, s'_i \leq s_{\sigma(i)} + s_{\sigma(i+q')} \mbox{ pour tout } i \in \unqpr \\
\quad \quad \bullet \, \, s''_\ell  = s_{\sigma(\ell+2q')} \mbox{ pour tout } \ell \in \unqsec.
\end{array}
\right.
$$
De plus, pour la combinaison lin\'eaire  construite dans la preuve   : 
\begin{itemize}
\item Les coefficients de la combinaison lin\'eaire peuvent \^etre calcul\'es de mani\`ere explicite 
et ils admettent $\dd_n ^{s_1 + \ldots + s_p}$ pour d\'enominateur commun.
\item Le coefficient du produit~\eqref{eq999} ne d\'epend que des $j_{\ell}$ et des $s_\ell$ 
pour $\ell \in \{ \sigma(2q' + q'' + 1), \ldots, \sigma(p)\}$.
\end{itemize}
\end{Th}

Dans ce th\'eor\`eme, et dans toute la suite, on identifie le groupe $\zdz$ \`a $\{-1, 1\}$ : pour 
$\eps \in \zdz$ et $s $ entier, on a $\eps^s = 1$ si $s$ est pair et $\eps^s = -1$ si $s$ est impair.

\bigskip

Les contraintes sur les produits \eqref{eq999} signifient que les polyz\^etas $\zetaetoile(s'_i)$ 
de profondeur 1   apparaissent par une sorte de concat\'enation de deux indices : c'est pourquoi ils 
peuvent appara\^{\i}tre jusqu'\`a $s'_i = 2A - 1$ dans les th\'eor\`emes \ref{thconj1} et \ref{thconj2}. 
C'est aussi la raison pour laquelle $q'$ appara\^{\i}t avec un facteur 2 dans la majoration $2q'+q'' \leq p$. 
En revanche, les $s''_\ell$ de \eqref{eq999} sont directement une sous-famille du $p$-uplet 
initial $(s_1, \ldots,  s_p)$ (\`a permutation pr\`es). La remarque qui termine l'\'enonc\'e du th\'eor\`eme 
\ref{thconj3} signifie que le coefficient de~\eqref{eq999} ne d\'epend ni des $s_\ell$ de cette sous-famille 
ni de ceux qui contr\^olent par concat\'enation les $s'_i$, mais seulement des autres (s'il y en a ; sinon, 
c'est que le coefficient ne d\'epend ni de $s_1, \ldots, s_p$ ni de $j_1, \ldots, j_p$). 

\bigskip

Si la profondeur $p$ est inf\'erieure ou \'egale \`a 3, les produits \eqref{eq999}  
sont des produits de valeurs de z\^eta en des entiers impairs,  ou bien des polyz\^etas 
antisym\'etriques de profondeur 2 ou 3. On va maintenant expliciter, \`a titre d'exemple,  
le coefficient d'un  tel polyz\^eta antisym\'etrique $\zetaetoileantisym (s''_{1}, \ldots, 
s''_{q''})$ dans la combinaison lin\'eaire \eqref{eq4}. La preuve de ce r\'esultat sera donn\'ee 
en m\^eme temps que celle du th\'eor\`eme \ref{thconj3}, aux paragraphes \ref{subsec52} et \ref{subsec53}.  

\smallskip

Si $p=2$, un tel polyz\^eta ne peut  appara\^{\i}tre (avec un coefficient non nul) que  si $s_1$ et $s_2$  sont impairs ; dans ce cas, sa contribution est toujours $4 ( \zeta(s_1, s_2) - \zeta(s_2, s_1))$.

Supposons maintenant que $p=3$. Alors  des polyz\^etas antisym\'etriques de profondeur 2 et 3 peuvent appara\^{\i}tre. En profondeur 3, la seule contribution possible est dans le cas o\`u  $s_1$, $s_2$ et $s_3$ sont impairs ; elle vaut 
$$8 \zetaetoileantisym (s_1, s_2, s_3) = 8 \sum_{\sigma \in \strois} \eps_\sigma \zetaetoile (s_{\sigma(1)},s_{\sigma(2)}, s_{\sigma(3)}).$$
Explicitons maintenant  la contribution des  polyz\^etas antisym\'etriques de  profondeur 2 (qui correspondent \`a $q'=0$ et $q'' = 2$). C'est une combinaison lin\'eaire des polyz\^etas $\zetaetoile (s_{i+1}, s_{i+2}) - \zetaetoile (s_{i+2}, s_{i+1})$ pour $i = 1, 2, 3$ (en interpr\'etant les indices modulo 3, par exemple $s_4 = s_1$). Ce polyz\^eta antisym\'etrique n'appara\^{\i}t que si $s_{i+1}$ et $s_{i+2}$ sont impairs. Dans ce cas, son coefficient est 
$$-4 \Big( \sum_{\ell = 1}^{j_i} \frac{1}{\ell^{s_i}} +  \sum_{\ell = 1}^{n-j_i} \frac{1}{\ell^{s_i}} \Big)$$
si $s_i$ est impair. Si $s_i$ est pair et $j_i \geq n/2$, c'est
$$-4 \Big( \sum_{\ell = n-j_i+1}^{j_i} \frac{1}{\ell^{s_i}}  \Big).$$
Enfin, si $s_i$ est pair et $j_i \leq n/2$, c'est
$$+4 \Big( \sum_{\ell = j_i+1}^{n-j_i} \frac{1}{\ell^{s_i}}  \Big).$$
Dans chacun de ces trois cas, on voit que ce coefficient ne d\'epend pas de $s_{i+1}$,  $s_{i+2}$,  
$j_{i+1}$,  $j_{i+2}$,  mais seulement de  $s_{i}$ et  de $j_{i}$ (comme \'enonc\'e dans le th\'eor\`eme \ref{thconj3}). 

\bigskip

\noindent{\bf Question : } Pourrait-on utiliser ces expressions explicites (en profondeur 2 ou 3)  
 pour trouver des polyn\^omes $P$ pour lesquels la partie ``polyz\^etas antisym\'etriques'' de la 
combinaison lin\'eaire du th\'eor\`eme \ref{thconj1} est nulle ? Pour ces polyn\^omes, cette 
combinaison lin\'eaire serait donc un polyn\^ome en valeurs de $\zeta$ en des entiers impairs. 

\subsection{Liens entre les th\'eor\`emes \ref{thconj2} et \ref{thconj3}} \label{subsec42bis}

Comme on va le voir, le th\'eor\`eme \ref{thconj3} est une forme plus pr\'ecise du th\'eor\`eme \ref{thconj2}.

\bigskip

Pour d\'eduire le  th\'eor\`eme  \ref{thconj2} du  th\'eor\`eme  \ref{thconj3}, on proc\`ede comme suit 
(il s'agit de la m\^eme strat\'egie que celle d\'etaill\'ee au paragraphe \ref{subsec43} ci-dessous). 
\'Etant donn\'e $P \in \mathscr{A}_p$, on utilise le d\'eveloppement en 
\'el\'ements simples du paragraphe \ref{subsec41} et on regroupe les termes 
qui correspondent \`a une m\^eme orbite sous l'action du groupe $\zdzp \psd \spp$ 
(voir \S \ref{subsec41}). Le fait que  $P \in \mathscr{A}_p$ signifie (voir 
\'egalement  \S \ref{subsec41}) que tous ces termes apparaissent 
avec le m\^eme coefficient, au signe pr\`es (et ce signe est donn\'e par la signature).
 On est donc ramen\'e \`a \'evaluer la somme sur chaque orbite, qui est exactement de la 
forme \eqref{eq4}: il suffit d'appliquer le th\'eor\`eme  \ref{thconj3}. 

\bigskip

R\'eciproquement, en mettant au m\^eme d\'enominateur les termes obtenus quand  $\sigma$ et  
$(\eps_1,\ldots,\eps_p)$ varient, on voit que \eqref{eq4} est de la forme \eqref{eq1} pour un 
certain polyn\^ome $P  \in \mathscr{A}_p$, de degr\'e $\leq A(n+1)-1$ par rapport \`a chacune 
des variables. Ceci prouve que le th\'eor\`eme  \ref{thconj2} implique le th\'eor\`eme \ref{thconj3}, 
\`a condition d'oublier,   dans ce dernier,  les pr\'ecisions donn\'ees en compl\'ement.

\subsection{Preuve que le th\'eor\`eme \ref{thconj3}  implique le th\'eor\`eme  \ref{thconj1}}
										\label{subsec43}

Commen\c cons par le point d\'elicat, qui diff\'erencie cette preuve de celle du paragraphe \ref{subsec42bis}. 

Sous les hypoth\`eses du th\'eor\`eme  \ref{thconj1}, la fraction rationnelle $R$ d\'efinie par 
\eqref{eq041} est de degr\'e $\leq -2$ par rapport \`a chacune de ses variables.  
Donc $k_1 R(k_1, \ldots, k_p)$ tend vers 0 quand $k_1$ tend vers l'infini, et on obtient en passant \`a la limite dans \eqref{eq2} :
$$\sum_\indso{0 \leq j_2, \ldots, j_p \leq  n}{1 \leq s_2, \ldots, s_p \leq A }
\frac{1}{(k_2+j_2)^{s_2} \ldots (k_p+j_p)^{s_p}} \sum_{j_1 = 0} ^n \cjsun = 0.$$
Par unicit\'e du d\'eveloppement en \'el\'ements simples de la fraction rationnelle nulle, 
on obtient pour tous $s_2, \ldots, s_p, j_2, \ldots, j_p $  :
$$ \sum_{j_1 = 0} ^n \cjsun = 0.$$
Le m\^eme raisonnement, appliqu\'e avec $k_i$ au lieu de $k_1$, montre que pour tout $i \in \unp$ on a :
 \begin{equation} \label{eq2pr}
 \sum_{j_i = 0} ^n \cjsinv = 0. 
 \end{equation}
 
 \bigskip
 
 Une fois ce r\'esultat pr\'eliminaire \'etabli, on peut suivre la strat\'egie r\'esum\'ee au 
paragraphe \ref{subsec42bis}, combin\'ee avec la r\'egularisation des divergences  logarithmiques 
et une \'etude plus d\'etaill\'ee de l'action du groupe (n\'ecessaire pour utiliser \eqref{eq2pr}).

\bigskip

En utilisant le d\'eveloppement en \'el\'ements simples \eqref{eq2}, on voit que la s\'erie 
convergente \eqref{eq021} est la limite, quand $N$ tend vers l'infini, de la somme
\begin{equation} \label{eq818}
\sum_\indso{0 \leq j_1, \ldots, j_p \leq  n}{1 \leq s_1, \ldots, s_p \leq A } \cjs
\sum_{N \geq k_1 \geq  \ldots \geq   k_p \ge 1} 
\frac{1}{(k_1+j_1)^{s_1} \ldots (k_p+j_p)^{s_p}}. 
 \end{equation}
Or l'ensemble d'indices $\zeron^p \croix \una^p$ est la r\'eunion disjointe des orbites sous l'action 
du groupe $\zdzp \psd \spp$ d\'efinie au paragraphe \ref{subsec41}. Etudions la contribution de chaque 
orbite \`a cette somme. Fixons  $(\jsoul, \ssoul) = (j_1, \ldots, j_p, s_1, \ldots, s_p)$, et consid\'erons 
un point quelconque  $(\jprsoul, \sprsoul) = (j'_1, \ldots, j'_p, s'_1, \ldots, s'_p)$ de son orbite 
(not\'ee $\orb$). Il existe $\gamma \in \spp$ et $(\eps_1, \ldots, \eps_p) \in \zdzp$ tels que 
$s'_1 = s_{\gamma(1)}$, \ldots, $s'_p = s_{\gamma(p)}$, $j'_1 = \eps_1 \cdot j_{\gamma(1)}$, \ldots,  
$j'_p = \eps_p \cdot j_{\gamma(p)}$. La relation \eqref{eq3} donne 
\begin{equation} \label{eq820}
\cjspr = \eps_\gamma  \eps_1^{s_{\gamma(1)}+1} \ldots \eps_p ^{s_{\gamma(p)}+1}  \cjs .
 \end{equation}

Comme tout \'el\'ement $(\jprsoul, \sprsoul) $ de $\orb$ s'\'ecrit ainsi pour exactement 
$\frac{2^p p!}{\Card \orb}$ \'el\'ements $(\eps_1, \ldots, \eps_p, \gamma) \in  \zdzp \psd  \spp$, 
on voit que la contribution de $\orb$ \`a la somme \eqref{eq818} est exactement la somme \eqref{eq4}, 
multipli\'ee par $\cjs \frac{\Card \orb}{2^p p!}$. D'apr\`es le th\'eor\`eme \ref{thconj3}, cette contribution s'\'ecrit donc 
\begin{equation} \label{eq819}
\frac{\Card \orb}{2^p p!} \cjs   Q_{\jsoul, \ssoul} (H_N) +   \gdo_\eps (N^{-1+\eps})
 \end{equation}
pour tout $\eps > 0$. Or pour $(\jprsoul, \sprsoul)  \in \orb$, en prenant $(\eps_1, \ldots, \eps_p, \gamma) $ comme ci-dessus, on voit par unicit\'e du polyn\^ome $Q_{\jprsoul, \sprsoul}$ que
$$Q_{\jprsoul, \sprsoul}(X) =  \eps_\gamma  \eps_1^{s_{\gamma(1)}+1} \ldots \eps_p ^{s_{\gamma(p)}+1}    Q_{\jsoul, \ssoul} (X). $$
Compte tenu de \eqref{eq820}, on peut donc \'ecrire \eqref{eq819} sous la forme
$$\sum_{(\jprsoul, \sprsoul) \in \orb} \frac{1}{2^p p!} \cjspr Q_{\jprsoul, \sprsoul} (H_N) + \gdo_\eps (N^{-1+\eps})$$
pour tout $\eps > 0$. Cette \'ecriture de la contribution de $\orb$ \`a \eqref{eq818} montre que la somme \eqref{eq818} est \'egale \`a 
\begin{equation} \label{eq817}
 \frac{1}{2^p p!} \sum_\indso{0 \leq j_1, \ldots, j_p \leq  n}{1 \leq s_1, \ldots, s_p \leq A } \cjs
 Q_{\jsoul, \ssoul} (H_N) + \gdo_\eps (N^{-1+\eps})
  \end{equation}
pour tout $\eps > 0$. Cette \'ecriture (qui consiste \`a r\'e\'ecrire \eqref{eq818} en moyennant sur chaque orbite, puis en appliquant le th\'eor\`eme \ref{thconj3} \`a chacune d'elles) est le point crucial qui va permettre maintenant de conclure, en appliquant la relation \eqref{eq2pr} d\'emontr\'ee au d\'ebut du paragraphe.
 
\bigskip
 
Comme la somme \eqref{eq817}  converge vers \eqref{eq021} quand $N$ tend vers l'infini, le polyn\^ome 
\begin{equation} \label{eq1001}
Q(X) = \frac{1}{2^p p!} \sum_\indso{0 \leq j_1, \ldots, j_p \leq  n}{1 \leq s_1, \ldots, s_p \leq A } \cjs
Q_{\jsoul, \ssoul} (X)
\end{equation}
est en fait constant, \'egal \`a sa valeur en 0; et cette valeur est exactement la somme \eqref{eq021}. Donc le th\'eor\`eme  \ref{thconj3} montre que \eqref{eq021} est une combinaison lin\'eaire, 
 \`a coefficients rationnels,  de produits de la forme
\begin{equation} \label{eq1000}
\zetaetoile(s'_1) \ldots  \zetaetoile(s'_{q'}) \zetaetoileantisym (s''_{1}, \ldots, s''_{q''})
\end{equation}
avec, pour chaque produit de cette forme,  
$q' ,q'' \geq 0$  entiers tels que $ 2q'+q'' \leq p$, 
$s'_1, \ldots, s'_{q'}, s''_{1}, \ldots, s''_{q''} \geq 1$    impairs, 
et $ \sigma \in \spp $  tel que $ s'_i \leq s_{\sigma(i)} + s_{\sigma(i+q')}$ pour tout $i \in \unqpr $ et $s''_\ell  = s_{\sigma(\ell+2q')} $ pour tout $  \ell \in \unqsec$.

\bigskip

Comme $\zetaetoile(1) = 0$, on peut supposer que dans un tel produit \eqref{eq1000} on a $s'_1, \ldots, s'_{q'} \geq 3$. Si on a aussi $s''_{1}, \ldots, s''_{q''} \geq 3$, alors ce produit fait partie de ceux autoris\'es dans la conclusion du th\'eor\`eme \ref{thconj1}, donc il n'y a rien d'autre \`a d\'emontrer. Supposons en revanche que $s''_\ell = s_{\sigma(\ell+2q')}=1$ pour un certain $\ell \in \unqsec$. D'apr\`es les pr\'ecisions donn\'ees \`a la fin du  th\'eor\`eme \ref{thconj3}, le coefficient du produit \eqref{eq1000} dans la d\'ecomposition de $Q_{\jsoul, \ssoul}(0)$    ne d\'epend pas de $j_{\sigma(\ell + 2q')}$. D'apr\`es \eqref{eq1001} et l'\'egalit\'e \eqref{eq2pr} d\'emontr\'ee au d\'ebut de ce paragraphe (appliqu\'ee avec $i = \sigma(\ell+2q')$), ce produit appara\^{\i}t dans $Q(0)$ avec un coefficient nul, donc ne contribue pas \`a la somme \eqref{eq021}.

\bigskip

Ceci termine la preuve du fait que le th\'eor\`eme \ref{thconj3}  implique le th\'eor\`eme  \ref{thconj1}.

\bigskip

Pour d\'emontrer l'assertion sur le d\'enominateur \label{subsecdenompreuve} des coefficients qui figure \`a la  fin du paragraphe~\ref{subsec22},  il suffit d'appliquer le lemme suivant et  de suivre, dans toute la preuve ci-dessus, les d\'enominateurs des nombres rationnels qui apparaissent.

\begin{Lemme} \label{lemdenom} Si $P = n!^{Ap} \widetilde P $ o\`u $\widetilde P $ est un polyn\^ome \`a coefficients entiers, alors 
$$\dd_n ^{Ap - (s_1+\ldots+s_p)} \cjs \in \Z$$
pour tous $s_1, \ldots, s_p, j_1, \ldots, j_p $.
\end{Lemme}

D\'emontrons maintenant ce lemme. Par $\Z$-lin\'earit\'e, il suffit de traiter le cas o\`u $P = n!^{Ap} X_1^{r_1} \ldots X_p^{r_p}$ avec $0 \leq r_1, \ldots, r_p \leq A(n+1)-1$. Admettons pour l'instant la propri\'et\'e suivante en une variable : pour tout $r \in \{0, \ldots, A(n+1)-1\}$ on a 
\begin{equation} \label{eqdenomunevar}
\frac{n! ^A k^r}{(k)_{n+1}^A} = \sum_{j=0}^n \sum_{s=1}^A \frac{\ejsr}{(k+j)^s} 
\mbox{ avec } \dd_n^{A-s} \ejsr \in \Z.
\end{equation}
 Le produit de cette relation, \'ecrite avec $k=k_i$ pour $i \in \{1, \ldots, p\}$, montre que la fraction rationnelle \eqref{eq041} peut s'\'ecrire  sous la forme \eqref{eq2} avec
$$\dd_n^{Ap-(s_1+\ldots+s_p)} \cjs = \prod_{i=1}^p \dd_n^{A-s_i} \ejsri \in \Z.$$
Ceci termine la preuve du lemme, en admettant la relation \eqref{eqdenomunevar}.

D\'emontrons maintenant cette relation. La matrice de passage de la base canonique $(k^r)_{0 \leq r \leq A(n+1)-1}$ \`a la base form\'ee par les polyn\^omes $(k)_{n+1}^a (k+n+1-\sigma)_{\sigma}$ (pour $0 \leq a \leq A-1$ et $0 \leq \sigma \leq n$) est \`a coefficients entiers, triangulaire sup\'erieure \`a diagonale de 1. Donc son inverse l'est aussi ; ceci permet de d\'ecomposer le mon\^ome $k^r$ dans la nouvelle base (avec des coefficients entiers). Par $\Z$-lin\'earit\'e, on est ramen\'e \`a d\'ecomposer des fractions rationnelles de la forme $\frac{n!^A}{(k)_{n+1}^{A-a-1} (k)_{n+1-\sigma}}$. Pour cela, on utilise la formule suivante   :
$$\frac{n!  }{(k)_{n+1}} = \sum_{j=0}^n  \frac{H_{n,j}}{k+j} \mbox{ avec } H_{n,j}  \in \Z.$$
Cette formule (qui est simplement le cas particulier $A=1$, $r=0$ de \eqref{eqdenomunevar}) est 
d\'emontr\'ee par exemple dans le lemme 5 de \cite{BR}. 
Il suffit alors de faire le produit cette formule, appliqu\'ee $a$ fois sous cette forme et une fois avec $n$ remplac\'e par $n-\sigma$. Une fois ce produit d\'evelopp\'e, on utilise (comme dans \cite{Colmez})  la formule $\frac{1}{(k+j)(k+j')} = \frac{1}{(j'-j)(k+j)} + \frac{1}{(j-j')(k+j')}$ pour $j \neq j'$. Chaque application de cette formule fait appara\^{\i}tre un d\'enominateur, qui est un diviseur de $\dd_n$. Apr\`es de multiples applications de cette formule, on arrive \`a une somme de la forme annonc\'ee dans   \eqref{eqdenomunevar}, et le coefficient $\ejsr$ est la somme de plusieurs termes qui proviennent tous d'au plus $A-s$ applications de cette formule. Ceci termine la preuve de \eqref{eqdenomunevar}, donc celle du lemme.

\section{D\'emonstration du th\'eor\`eme \ref{thconj3}  } \label{sec5}

Dans cette partie, on d\'emontre le  th\'eor\`eme \ref{thconj3}  par r\'ecurrence sur la profondeur. En th\'eorie, l'initialisation (\S \ref{subsec51}) et le c\oe ur de la r\'ecurrence (\S \ref{subsec54}) suffisent ; mais on d\'emontre aussi compl\`etement les cas $p=2$ (\S \ref{subsec52}) et $p=3$ (\S \ref{subsec53}) pour illustrer et motiver les constructions du paragraphe \ref{subsec54}.

C'est dans cette partie, et nulle part ailleurs, que l'hypoth\`ese ``$n$ est pair'' est utilis\'ee (voir la remarque \ref{rqparite} ci-dessous). 

\subsection{Preuve du th\'eor\`eme \ref{thconj3}  en profondeur 1} \label{subsec51}

Quand $p=1$, le th\'eor\`eme \ref{thconj3} concerne des s\'eries de la forme 
$$\sum_{k=1} ^N \Big(  \frac{1}{(k+j)^s} + \frac{(-1)^{s+1}}{(k+n-j)^s} \Big).$$
Si $s \geq 2$, on voit directement que cette somme vaut $(1+(-1)^{s+1}) \zeta(s) + \rho_{j,s} + \gdo(\frac{1}{N})$ avec $\dd_n^s \rho_{j,s} \in \Z$, ce qui d\'emontre le th\'eor\`eme dans ce cas. Sinon, c'est-\`a-dire si $s=1$, cette somme vaut $2 H_N  +   \rho_{j,s} + \gdo(\frac{1}{N})$ avec $\dd_n \rho_{j,s} \in \Z$, ce qui d\'emontre aussi le r\'esultat voulu puisque $\zetaetoile(1)=0$.

Le th\'eor\`eme \ref{thconj3}  est donc d\'emontr\'e quand $p=1$.

\subsection{Preuve du th\'eor\`eme \ref{thconj3}  en profondeur 2} \label{subsec52}

Dans ce paragraphe, on suppose $p=2$ et on d\'emontre, par r\'ecurrence sur $(j_1,j_2)$, que le th\'eor\`eme \ref{thconj3}  est vrai pour tous $s_1, s_2$. L'entier  $n$ est fix\'e dans toute la preuve. 

\smallskip

L'initialisation de cette r\'ecurrence est le cas o\`u $j_1 = j_2 = \frac{n}{2}$ (puisque $n$ est suppos\'e pair ; voir la remarque \ref{rqparite} ci-dessous).  La somme \eqref{eq4} est alors nulle si $s_1 $ ou $s_2$ est pair ; le r\'esultat du th\'eor\`eme \ref{thconj3} est trivial dans cette situation. On peut donc supposer que $s_1$ et $s_2$ sont impairs. La somme \eqref{eq4} vaut alors $4(\tau_{s_1, s_2} - \tau_{s_2, s_1})$, en posant 
$$\tau_{s_1, s_2} = \sum_{N \geq k_1 \geq k_2 \geq 1} \frac{1}{(k_1 + \frac{n}{2})^{s_1}  (k_2 + \frac{n}{2})^{s_2}}.$$
Or on a 
\begin{eqnarray*}
\tau_{s_1, s_2} 
&=&  \sum_{N + \frac{n}{2} \geq \ell_1 \geq \ell_2 \geq  \frac{n}{2} + 1} \frac{1}{\ell_1^{s_1} \ell_2^{s_2}}\\
&=&  \sum_{N + \frac{n}{2} \geq \ell_1 > \ell_2 \geq  1} \frac{1}{\ell_1^{s_1} \ell_2^{s_2}} 
+ \sum_{N + \frac{n}{2} \geq \ell  \geq  1} \frac{1}{\ell^{s_1+s_2}} 
- \sum_{\ell_2 = 1} ^{n/2}  
\frac{1}{\ell_2^{s_2}} \sum_{\ell_1 = \ell_2} ^{N + \frac{n}{2}} \frac{1}{\ell_1^{s_1}}.
\end{eqnarray*}
La proposition \ref{propstuffle} fournit un polyn\^ome $Q$ tel que, pour tout $\eps > 0$, on ait quand $N$ tend vers l'infini : 
$$\tau_{s_1, s_2}  = Q(H_N) +  \gdo_\eps (N^{-1+\eps})$$
avec 
$$Q(0) = \zetaetoile (s_1,s_2) + \zeta(s_1+s_2) - \Big( \sum_{\ell_2 = 1} ^{n/2} \frac{1}{\ell_2 ^{s_2}} \Big)  \zetaetoile(s_1) + r ,$$
o\`u $\dd_n ^{s_1+s_2} r \in \Z$ (ceci provient du fait que $ Q(H_{N+n/2}) = Q(H_N) +  \gdo_\eps (N^{-1+\eps})$).
Or $\zetaetoile(1) = 0$ et $\zetaetoile(s_1) = \zeta(s_1)$ pour $s_1 \geq 2$. \'donn\'e que la somme \eqref{eq4} vaut  $4(\tau_{s_1, s_2} - \tau_{s_2, s_1})$, cela d\'emontre le th\'eor\`eme \ref{thconj3}  quand $s_1$ et $s_2$ sont impairs, avec $j_1 = j_2 = n/2$. Cela termine la preuve de l'initialisation de la r\'ecurrence.

\bigskip

\begin{remark} \label{rqparite} Dans cette initialisation, on a suppos\'e que $n$ est pair. C'est le seul endroit dans cet article (avec les initialisations analogues en profondeurs $3$ et $p \geq 4$ aux paragraphes \ref{subsec53} et \ref{subsec54})  o\`u cette hypoth\`ese est utilis\'ee. Si on voulait d\'emontrer les m\^emes r\'esultats lorsque $n$ est impair, il suffirait de d\'emontrer cette initialisation dans ce cas. Bien entendu, on ne pourrait plus prendre $j_1 = j_2 = \frac{n}{2}$, donc les calculs seraient plus compliqu\'es. On pourrait par exemple choisir  $j_1 = j_2 =0$. 
\end{remark}

\bigskip
 
 La suite de la d\'emonstration consiste \`a \'etablir le r\'esultat suivant pour tous $j_1  \in \zeronmu$ et $j_2 \in \zeron$  :
\begin{equation}  \label{eqrec2}
\left\{
\begin{array}{l}
\mbox{ le th\'eor\`eme \ref{thconj3}   est vrai pour le couple $(j_1, j_2)$, quels que soient  $s_1$ et $ s_2$}, \\
\quad \quad \mbox{ si, et seulement si, }\\
\mbox{il  est vrai pour le couple $(j_1+1, j_2)$, quels que soient  $s_1$ et $ s_2$}.
\end{array}
\right.
\end{equation}
En effet, supposons \eqref{eqrec2} \'etablie. Comme  le th\'eor\`eme  est vrai pour le couple $(j_1 = n/2, j_2 = n/2)$, il est vrai pour $(j_1, n/2)$ quel que soit $j_1 \in \unn$ en utilisant \eqref{eqrec2}. Or quand on \'echange $j_1$ et $j_2$, ainsi que (simultan\'ement)$s_1$ et $s_2$, la somme \eqref{eq4} est chang\'ee en son oppos\'e. Donc  le th\'eor\`eme  est vrai pour $(j_1, j_2)$ et $(s_1, s_2)$ si, et seulement si, il est vrai pour $(j_2, j_1)$ et $(s_2, s_1)$. En particulier,  le th\'eor\`eme  est donc vrai 
pour $(n/2, j_2)$ quel que soit $j_2 \in \unn$, et quels que soient $s_1$ et $s_2$. En appliquant \`a nouveau \eqref{eqrec2}, on voit que  le th\'eor\`eme  est vrai pour tout couple $(j_1, j_2) \in \unn ^2$.

Pour terminer la preuve du th\'eor\`eme \ref{thconj3}   en profondeur 2, il suffit donc d'\'etablir \eqref{eqrec2}. 

\bigskip

Posons 
$$K_N(j_1, j_2, s_1, s_2) = \sum_{N \geq k_1 \geq k_2 \geq 1}  \frac{1}{(k_1+j_1)^{s_1} (k_2+j_2)^{s_2}}.$$
Alors le th\'eor\`eme \ref{thconj3} concerne la somme
\begin{equation} \label{eq6}
\sum_{\eps_1, \eps_2 \in \zdz} \eps_1 ^{s_1+1}  \eps_2 ^{s_2+1} 
\Big( K_N(\eps_1 \cdot j_1, \eps_2 \cdot  j_2, s_1, s_2) -
K_N(\eps_2 \cdot j_2, \eps_1 \cdot  j_1, s_2, s_1)  \Big) .
\end{equation}
Pour \'etablir \eqref{eqrec2}, il suffit de d\'emontrer que la diff\'erence entre \eqref{eq6} pour $(j_1+1, j_2)$ et  \eqref{eq6} pour $(j_1, j_2)$ est de la forme annonc\'ee dans le th\'eor\`eme \ref{thconj3}. Pour \'evaluer cette diff\'erence, on aura besoin des calculs suivants.

D'abord, 
\begin{eqnarray}
\lefteqn{K_N(j'_1+1 , j'_2, s_1, s_2) - K_N(j'_1, j'_2, s_1, s_2)}\qquad \nonumber \\
&=& \sum_{N \geq k_1 \geq k_2 \geq 1} \frac{1}{(k_2+j'_2)^{s_2}} \Big(  \frac{1}{(k_1+j'_1+1)^{s_1}} -  \frac{1}{(k_1+j'_1)^{s_1}} \Big) \nonumber\\
&=& \sum_{k_2= 1} ^N  \frac{1}{(k_2+j'_2)^{s_2}}  \sum_{k_1 = k_2} ^N 
\Big(  \frac{1}{(k_1+j'_1+1)^{s_1}} -  \frac{1}{(k_1+j'_1)^{s_1}} \Big) \nonumber\\
&=& \sum_{k_2= 1} ^N  \frac{1}{(k_2+j'_2)^{s_2}}   
\Big(  \frac{1}{(N+j'_1+1)^{s_1}} -  \frac{1}{(k_2+j'_1)^{s_1}} \Big) \nonumber\\
&=& \gdo (\frac{\log N}{N}) - \sum_{k= 1} ^N  \frac{1}{(k+j'_1)^{s_1} (k+j'_2)^{s_2}}.   \label{eq7}
\end{eqnarray} 
On peut en d\'eduire, ou bien d\'emontrer de mani\`ere analogue, la relation
\begin{multline}
K_N(j'_1-1 , j'_2, s_1, s_2) - K_N(j'_1, j'_2, s_1, s_2)
\\
= \sum_{k_2= 1} ^N  \frac{1}{(k_2+j'_2)^{s_2}}   
\Big(  \frac{-1}{(N+j'_1)^{s_1}} +  \frac{1}{(k_2+j'_1-1)^{s_1}} \Big) 
\\
= \gdo (\frac{\log N}{N}) + \sum_{k= 1} ^N  \frac{1}{(k+j'_1-1)^{s_1} (k+j'_2)^{s_2}}.  \label{eq8}
\end{multline} 
On aura aussi besoin des relations suivantes, dont la preuve est analogue, et dans lesquelles c'est la deuxi\`eme variable que l'on modifie : 
\begin{multline}
K_N(j'_2, j'_1+1 , s_2, s_1) - K_N(j'_2, j'_1, s_2, s_1)   \\
= \sum_{k= 1} ^N  \frac{1}{(k+j'_1+1)^{s_1} (k+j'_2)^{s_2}} 
- \frac{1}{(j'_1+1)^{s_1}} \sum_{k=1}^N \frac{1}{(k+j'_2)^{s_2}}   \label{eq12}
\end{multline} 
et
\begin{multline}
K_N(j'_2, j'_1-1 , s_2, s_1) - K_N(j'_2, j'_1, s_2, s_1)\\
= - \sum_{k= 1} ^N  \frac{1}{(k+j'_1)^{s_1} (k+j'_2)^{s_2}} 
+ \frac{1}{{j'_1}^{s_1}} \sum_{k=1}^N \frac{1}{(k+j'_2)^{s_2}} .   \label{eq13}
\end{multline} 

\bigskip

Posons 
$$\Delta_{\eps_1, \eps_2} (j_1, j_2) = 
K_N(\eps_1 \cdot (j_1+1) , \eps_2 \cdot  j_2, s_1, s_2)
- K_N(\eps_1 \cdot j_1, \eps_2 \cdot  j_2, s_1, s_2)$$
et 
$$ \Deltati_{\eps_1, \eps_2} (j_1, j_2) = 
K_N(\eps_2 \cdot j_2 , \eps_1 \cdot ( j_1+1), s_2, s_1)
-K_N(\eps_2 \cdot j_2 , \eps_1 \cdot  j_1, s_2, s_1). $$
Avec ces notations, la diff\'erence entre \eqref{eq6} pour $(j_1+1, j_2)$ et  \eqref{eq6} pour $(j_1, j_2)$ (que l'on cherche \`a \'evaluer) est 
\begin{equation} \label{eq11}
\sum_{\eps_1, \eps_2 \in \zdz} \eps_1 ^{s_1+1}  \eps_2 ^{s_2+1} 
\Big(  \Delta_{\eps_1, \eps_2} (j_1, j_2) 
-  \Deltati_{\eps_1, \eps_2} (j_1, j_2)    \Big) .
\end{equation}
Or on a :
$$\eps_1 \cdot ( j_1+1) = 
\left\{
\begin{array}{l}
(\eps_1 \cdot j_1) + 1 \mbox{ si } \eps_1 = +1 \\
(\eps_1 \cdot j_1) - 1 \mbox{ si } \eps_1 = -1.  
\end{array}
\right.
$$
En utilisant successivement deux fois \eqref{eq7}, deux fois \eqref{eq8}, deux fois \eqref{eq12} et deux fois \eqref{eq13}, on voit que \eqref{eq11} est la somme des huit termes suivants : 
\begin{align}
\Delta_{+1, +1 } (j_1, j_2) 
&= - \sum_{k= 1} ^N  \frac{1}{(k+j_1)^{s_1} (k+j_2)^{s_2}}  +  \gdo (\frac{\log N}{N}) , 
														\label{eq14}\\
(-1)^{s_2+1} \Delta_{+1, -1 } (j_1, j_2) 
&= (-1)^{s_2} \sum_{k= 1} ^N  \frac{1}{(k+j_1)^{s_1} (k+n-j_2)^{s_2}}  +  \gdo (\frac{\log N}{N}) ,
														\label{eq15}\\
(-1)^{s_1+1} \Delta_{-1, +1 } (j_1, j_2) 
&= (-1)^{s_1+1} \sum_{k= 1} ^N  \frac{1}{(k+n-j_1-1)^{s_1} (k+j_2)^{s_2}}  +  \gdo (\frac{\log N}{N}) ,
														\label{eq16}\\
(-1)^{s_1+s_2} \Delta_{-1, -1 } (j_1, j_2) 
&= (-1)^{s_1+s_2} \sum_{k= 1} ^N  \frac{1}{(k+n-j_1-1)^{s_1} (k+n-j_2)^{s_2}} \label{eq17} 
\\
& \hspace{7.5cm}  + 
\gdo (\frac{\log N}{N}),  \nonumber
\end{align}													
\begin{equation}													
- \Deltati_{+1, +1 } (j_1, j_2)  
= - \sum_{k= 1} ^N  \frac{1}{(k+j_1+1)^{s_1} (k+j_2)^{s_2}}  +  
\frac{1}{(j_1+1)^{s_1}} \sum_{k=1} ^N \frac{1}{(k+j_2)^{s_2}},	\label{eq18}
\end{equation}	
\begin{multline}													
(-1)^{s_2} \Deltati_{+1, -1 } (j_1, j_2) 
\\
= (-1)^{s_2} \sum_{k= 1} ^N  \frac{1}{(k+j_1+1)^{s_1} (k+n-j_2)^{s_2}}  +  
\frac{(-1)^{s_2+1}}{(j_1+1)^{s_1}} \sum_{k=1} ^N \frac{1}{(k+n-j_2)^{s_2}},	
\label{eq19}
\end{multline}
\begin{multline}
(-1)^{s_1} \Deltati_{-1, +1 } (j_1, j_2) 
\\
= (-1)^{s_1+1} \sum_{k= 1} ^N  \frac{1}{(k+n-j_1)^{s_1} (k+j_2)^{s_2}}  +  
\frac{(-1)^{s_1}}{(n-j_1)^{s_1}} \sum_{k=1} ^N \frac{1}{(k+j_2)^{s_2}},	\label{eq20}
\end{multline}
\begin{multline}
(-1)^{s_1+s_2+1} \Deltati_{-1, -1 } (j_1, j_2) 
\\
= (-1)^{s_1+s_2} \sum_{k= 1} ^N  \frac{1}{(k+n-j_1)^{s_1} (k+n-j_2)^{s_2}}  + 
\frac{(-1)^{s_1+s_2+1}}{(n-j_1)^{s_1}} \sum_{k=1} ^N \frac{1}{(k+n-j_2)^{s_2}}. 	\label{eq21}
\end{multline}	
												
On va montrer que la somme de ces huit quantit\'es est bien de la forme voulue, c'est-\`a-dire s'\'ecrit $Q(H_N) + \gdo_{\eps} (N^{-1+\eps})$ pour un certain polyn\^ome $Q$ dont la valeur en 0 est une combinaison lin\'eaire des polyz\^etas autoris\'es. Pour cela, on groupe les termes de la mani\`ere suivante : 
\begin{enumerate}
\item Le premier terme de  \eqref{eq18} avec \eqref{eq17}.
\item Le premier terme de  \eqref{eq19} avec \eqref{eq16}.
\item Le premier terme de  \eqref{eq20} avec \eqref{eq15}.
\item Le premier terme de  \eqref{eq21} avec \eqref{eq14}.
\item Le second terme de  \eqref{eq18} avec celui de \eqref{eq19}.
\item Le second terme de  \eqref{eq20} avec celui de \eqref{eq21}.
\end{enumerate}
Pour chacun de ces six groupements, il suffit d'appliquer le th\'eor\`eme \ref{thconj2} en profondeur 1 (c'est-\`a-dire essentiellement le th\'eor\`eme \ref{thwp}, qui est le ph\'enom\`ene de sym\'etrie habituel : voir \S \ref{subsec51}) pour conclure. 

\bigskip

Ceci termine la preuve de le th\'eor\`eme \ref{thconj3} en profondeur 2.

\subsection{Preuve du th\'eor\`eme \ref{thconj3}  en profondeur 3} \label{subsec53}

On proc\`ede par r\'ecurrence, comme au paragraphe \ref{subsec52}.

Pour initialiser la r\'ecurrence, on consid\`ere (puisque $n$ est suppos\'e pair, voir la remarque \ref{rqparite}) le cas o\`u $j_1 = j_2 = j_3 = n/2$. Dans ce cas, \eqref{eq4} vaut 0 si l'un au moins des $s_i$ est pair. Il ne reste donc \`a traiter que le cas o\`u les trois $s_i$ sont impairs. Dans ce cas, on a :
$$\eqref{eq4} = 8 \sum_{\sigma \in \strois} \eps_\sigma \tau_{s_{\sigma(1)},s_{\sigma(2)}, s_{\sigma(3)}}$$
avec
$$\tau_{s_1, s_2,s_3} = \sum_{N \geq k_1 \geq k_2 \geq k_3 \geq 1} \frac{1}{(k_1 + \frac{n}{2})^{s_1}  (k_2 + \frac{n}{2})^{s_2}  (k_3 + \frac{n}{2})^{s_3}   }.$$
Or on a  
\begin{eqnarray*}
\tau_{s_1, s_2, s_3} 
&=&  \sum_{N + \frac{n}{2} \geq \ell_1 \geq \ell_2  \geq \ell_3 \geq  \frac{n}{2} + 1} \frac{1}{\ell_1^{s_1} \ell_2^{s_2}  \ell_3^{s_3}   }\\
&=&  \sum_{N + \frac{n}{2} \geq \ell_1 \geq \ell_2  \geq \ell_3 \geq   1} \frac{1}{\ell_1^{s_1} \ell_2^{s_2}  \ell_3^{s_3}   } -   \sum_{\ell_3 = 1} ^{n/2}  \frac{1}{\ell_3^{s_3}}  
 \sum_{N + \frac{n}{2} \geq \ell_1 \geq \ell_2 \geq  \ell_3} \frac{1}{\ell_1^{s_1} \ell_2^{s_2}}  \\
&=&  \sum_{N + \frac{n}{2} \geq \ell_1 \geq \ell_2  \geq \ell_3 \geq   1} \frac{1}{\ell_1^{s_1} \ell_2^{s_2}  \ell_3^{s_3}   } -  \Big(  \sum_{\ell_3 = 1} ^{n/2}  \frac{1}{\ell_3^{s_3}}   \Big) 
 \sum_{N + \frac{n}{2} \geq \ell_1 \geq \ell_2 \geq  1} \frac{1}{\ell_1^{s_1} \ell_2^{s_2}} \\
&& \quad +  \Big( \sum_{\frac{n}{2} \geq \ell_3 > \ell_2 \geq  1} \frac{1}{\ell_3^{s_3} \ell_2^{s_2}} \Big)
 \sum_{\ell_1 = 1} ^{N + n/2}  \frac{1}{\ell_1^{s_1}}
 -   \sum_{ \frac{n}{2} \geq  \ell_3 >  \ell_2 >  \ell_1 \geq  1} \frac{1}{\ell_3^{s_3} \ell_2^{s_2}  \ell_1^{s_1}   }
\end{eqnarray*} 
donc $ \tau_{s_1, s_2, s_3}  = Q(H_N) + \gdoeps$ pour un certain polyn\^ome $Q$ tel que (d'apr\`es la proposition \ref{propstuffle}) : 
\begin{eqnarray*}
Q(0) &=& \zetaetoile(s_1,s_2,s_3) + \zeta(s_1+s_2, s_3) + \zetaetoile(s_1, s_2+s_3) + \zeta(s_1+s_2+s_3)\\
&& - \Big( \sum_{\ell_3 = 1} ^{n/2} \frac{1}{\ell_3^{s_3}} \Big) 
\Big( \zetaetoile (s_1, s_2) +  \zeta (s_1+s_2) \Big)  + \chi(s_3, s_2) \zetaetoile(s_1) 
-   \sum_{ \frac{n}{2} \geq  \ell_3 >  \ell_2 >  \ell_1 \geq  1} \frac{1}{\ell_3^{s_3} \ell_2^{s_2}  \ell_1^{s_1} } 
\end{eqnarray*}
en posant
$$\chi(s_3, s_2)  = \sum_{\frac{n}{2} \geq \ell_3 > \ell_2 \geq  1} \frac{1}{\ell_3^{s_3} \ell_2^{s_2}}.$$
Ainsi, on obtient que \eqref{eq4} s'\'ecrit sous la forme $Q_1(H_N) + \gdoeps$ pour un certain polyn\^ome $Q_1$ tel que 
\begin{multline*}
Q_1(0) 
= 8 \sum_{\sigma \in \strois} \eps_\sigma \zetaetoile (s_{\sigma(1)},s_{\sigma(2)}, s_{\sigma(3)}) 
 - 8 \Big( \sum_{\ell = 1} ^{n/2} \frac{1}{\ell^{s_3}} \Big) \Big( \zetaetoile (s_1, s_2) 
 -  \zetaetoile (s_2, s_1) \Big) 
\\
 + 8 \Big( \sum_{\ell = 1} ^{n/2} \frac{1}{\ell^{s_2}} \Big) \Big( \zetaetoile (s_1, s_3)  -  \zetaetoile (s_3, s_1) \Big) 
 - 8 \Big( \sum_{\ell = 1} ^{n/2} \frac{1}{\ell^{s_1}} \Big) \Big( \zetaetoile (s_2, s_3)  -  \zetaetoile (s_3, s_2) \Big) 
\\
+ 8 \Big( \chi(s_3, s_2) - \chi(s_2, s_3)  \Big)  \zetaetoile(s_1)  
 - 8 \Big( \chi(s_3, s_1) - \chi(s_1, s_3)  \Big)  \zetaetoile(s_2)  
\\
- 8 \Big( \chi(s_1, s_2) - \chi(s_2, s_1)  \Big)  \zetaetoile(s_3)  
- 8  \sum_{\sigma \in \strois} \eps_\sigma 
 \sum_{ \frac{n}{2} \geq  \ell_3 >  \ell_2 >  \ell_1 \geq  1} \frac{1}{\ell_3^{s_{\sigma(3)}} \ell_2^{s_{\sigma(2)}}  \ell_1^{s_{\sigma(1)}}}. 
\end{multline*}
 Ceci termine l'initialisation de la r\'ecurrence.
 
 \bigskip 

D\'emontrons maintenant l'h\'er\'edit\'e. Pour raccourcir les notations, on pose $\jsoul = (j_1, j_2, j_3)$, $\ssoul = (s_1, s_2, s_3)$ et $\epssoul = (\eps_1, \eps_2, \eps_3)$. La preuve est parall\`ele \`a celle dans le cas de la profondeur 2 (\S  \ref{subsec52}), mais le groupement des termes qui permet de conclure est plus compliqu\'e. 

On pose 
$$K_N(\jsoul, \ssoul) = \sum_{N \geq k_1 \geq k_2  \geq k_3  \geq 1}  \frac{1}{(k_1+j_1)^{s_1} (k_2+j_2)^{s_2} (k_3+j_3)^{s_3}}$$
puis, pour $\sigma \in \strois$ : 
$$K_N^\sigma (\jsoul, \ssoul)  = K_N (j_{\sigma(1)}, j_{\sigma(2)}, j_{\sigma(3)}, s_{\sigma(1)}, s_{\sigma(2)}, s_{\sigma(3)})$$
de telle sorte que $K_N^{\Id} (\jsoul, \ssoul) = K_N  (\jsoul, \ssoul) $. Puisque $\eps_{\sigma^{-1}} = \eps_{\sigma}$, on a : 
$$\eqref{eq4} = \sum_{\epssoul \in \zdzt} \epssoul ^{\ssoul + 1} \sum_{\sigma \in \strois} \eps_\sigma
K_N^\sigma (\eps_1 \cdot j_1, \eps_2 \cdot j_2, \eps_3 \cdot j_3, \ssoul)$$
o\`u on note 
$\epssoul ^{\ssoul + 1}  = \eps_1 ^{s_1+ 1} \eps_2 ^{s_2+ 1} \eps_3 ^{s_3+ 1} .$
On pose aussi 
$$\Delta_{\epssoul} ^\sigma (\jsoul) = 
K_N^\sigma (\eps_1 \cdot (j_1+1), \eps_2 \cdot j_2, \eps_3 \cdot j_3, \ssoul)
- K_N^\sigma (\eps_1 \cdot j_1, \eps_2 \cdot j_2, \eps_3 \cdot j_3, \ssoul).$$
Alors la diff\'erence entre \eqref{eq4} pour $(j_1+1, j_2, j_3)$ et \eqref{eq4} pour $(j_1, j_2, j_3)$ est :
\begin{equation} \label{eq12nv} 
\sum_{\epssoul \in \zdzt} \epssoul ^{\ssoul + 1} \sum_{\sigma \in \strois} \eps_\sigma \Delta_{\epssoul} ^\sigma (\jsoul) .
\end{equation}
La suite de la preuve est consacr\'ee \`a \eqref{eq12nv} : il s'agit de montrer que cette somme est de la forme voulue, ce qui terminera la r\'ecurrence (de mani\`ere analogue \`a \eqref{eqrec2} dans le cas de la profondeur 2). Cette somme comprend 48 termes. Dans un premier temps, on fixe $\epssoul \in \zdzt$ et on explicite les 6 termes correspondants. Pour cela, on pose $j'_1 = \eps_1 \cdot j_1$, $j'_2 = \eps_2 \cdot j_2$,  $j'_3 = \eps_3 \cdot j_3$. Supposons d'abord que $\eps_1 =  + 1$ ; on a dans ce cas $\eps_1 \cdot (j_1+1) = j'_1 + 1$, et les six termes qui apparaissent correspondent aux formules  \eqref{eq7} et \eqref{eq12} du \S \ref{subsec52}. 

Commen\c cons par le terme qui provient du 3-cycle (123), qui envoie 1 sur 2, 2 sur 3 et 3 sur 1 : 
\begin{eqnarray}
\lefteqn{\Delta_{\epssoul} ^{(123)} (\jsoul)}\nonumber \qquad \\
&=& K_N(j'_2 , j'_3, j'_1+1,  s_2, s_3, s_1) - K_N(j'_2 , j'_3, j'_1,  s_2, s_3, s_1) \nonumber \\
&=& \sum_{N \geq k_1 \geq k_2 \geq k_3 \geq 1} \frac{1}{(k_1+j'_2)^{s_2}(k_2+j'_3)^{s_3}} \Big(  \frac{1}{(k_3+j'_1+1)^{s_1}} -  \frac{1}{(k_3+j'_1)^{s_1}} \Big) \nonumber\\
&=& \sum_{N \geq k_1 \geq k_2  \geq 1} \frac{1}{(k_1+j'_2)^{s_2}(k_2+j'_3)^{s_3}} \Big(  
\frac{1}{(k_2+j'_1+1)^{s_1}} -  \frac{1}{(j'_1+1)^{s_1}} \Big) \nonumber\\
&=& \sum_{N \geq k \geq \ell  \geq 1} \frac{1}{(\ell+j'_1+1)^{s_1}(k+j'_2)^{s_2}(\ell+j'_3)^{s_3}}  \nonumber\\
&& \hspace{4cm} - \frac{1}{(j'_1+1)^{s_1}}  \sum_{N \geq k \geq \ell  \geq 1} \frac{1}{(k+j'_2)^{s_2}(\ell+j'_3)^{s_3}} .
															\label{eq30}
\end{eqnarray} 
Ce terme appara\^{\i}t dans la somme \eqref{eq12nv} avec le coefficient $\epssoul^{\ssoul+1}$ (sous l'hypoth\`ese que $\eps_1  = +1$), de m\^eme que les cinq termes suivants, qui se calculent de mani\`ere analogue : 
\begin{align}   
 \Delta_{\epssoul} ^{\Id} (\jsoul)
&= -  \sum_{N \geq k \geq \ell  \geq 1} \frac{1}{(k+j'_1)^{s_1}(k+j'_2)^{s_2}(\ell+j'_3)^{s_3}} 
+ \gdo ( \frac{\log^2 N}{N})\label{eq34}  \\
- \Delta_{\epssoul} ^{(23)} (\jsoul)
&=  \sum_{N \geq k \geq \ell  \geq 1} \frac{1}{(k+j'_1)^{s_1}(\ell +j'_2)^{s_2}(k+j'_3)^{s_3}} 
+ \gdo ( \frac{\log^2 N}{N})  \label{eq33}  \\
- \Delta_{\epssoul} ^{(13)} (\jsoul) 
&= -  \sum_{N \geq k \geq \ell  \geq 1} \frac{1}{(\ell+j'_1+1)^{s_1}(\ell+j'_2)^{s_2}(k+j'_3)^{s_3}} \label{eq32} \\
& \qquad \qquad + \frac{1}{(j'_1+1)^{s_1}}  \sum_{N \geq k \geq \ell  \geq 1} \frac{1}{(\ell+j'_2)^{s_2}(k+j'_3)^{s_3}}	\nonumber
\end{align}													
\begin{align}
- \Delta_{\epssoul} ^{(12)} (\jsoul) 
&= -  \sum_{N \geq k \geq \ell  \geq 1} \frac{1}{(k+j'_1+1)^{s_1}(k+j'_2)^{s_2}(\ell+j'_3)^{s_3}} \label{eq31} \\
& \qquad \qquad +    \sum_{N \geq k \geq \ell  \geq 1} \frac{1}{(\ell+j'_1)^{s_1}(k+j'_2)^{s_2}(\ell+j'_3)^{s}}\nonumber
\\	
\Delta_{\epssoul} ^{(132)} (\jsoul) 
&=   \sum_{N \geq k \geq \ell  \geq 1} \frac{1}{(k+j'_1+1)^{s_1}(\ell+j'_2)^{s_2}(k+j'_3)^{s_3}} 	\label{eq35} \\
&  \qquad \qquad -    \sum_{N \geq k \geq \ell  \geq 1} \frac{1}{(\ell+j'_1)^{s_1}(\ell+j'_2)^{s_2}(k+j'_3)^{s_3}} .	\nonumber
\end{align}															
Si $\eps_1=-1$, il suffit de prendre l'oppos\'e du membre de droite, et d'y remplacer $j'_1$ par $j'_1-1$, pour que les formules \eqref{eq30} \`a \eqref{eq35} soient correctes. Les formules ainsi obtenues sont les analogues de \eqref{eq8} et \eqref{eq13} (au \S \ref{subsec52}). 
Pour ne pas  avoir \`a distinguer suivant la valeur de $\eps_1$, on aurait pu multiplier le membre de droite par $\eps_1$, et y remplacer $j'_1$ par $\jprun$. Gr\^ace \`a ces modifications, les formules 
 \eqref{eq30} \`a \eqref{eq35} auraient \'et\'e valables quel que soit $\eps_1$ ; on utilisera cette convention dans la suite.
 
Pour exprimer \eqref{eq12nv} sous une forme exploitable, on groupe deux par deux les termes obtenus, par les formules  \eqref{eq30} \`a \eqref{eq35}, \`a partir des 48 termes de la somme \eqref{eq12nv}. Comme \eqref{eq30} et \eqref{eq33} ne donnent qu'un terme (\`a part le terme d'erreur, qu'on omet dans toute la suite des calculs), et que \eqref{eq34}, \eqref{eq32}, \eqref{eq31} et \eqref{eq35} en donnent deux, on \'ecrit ainsi \eqref{eq12nv} comme une somme de $8 \times 10 = 80$ termes. On va maintenant expliciter ces 40 groupes de 2 termes.

Soit $\epssoul \in \zdzt$. On pose comme ci-dessus $j'_1 = \eps_1 \cdot j_1$, $j'_2 = \eps_2 \cdot j_2$,  $j'_3 = \eps_3 \cdot j_3$.  Les 5 groupes qui correspondent \`a $\epssoul$ sont les suivants :
\begin{enumerate}
\item On regroupe le deuxi\`eme terme de \eqref{eq30} avec celui de \eqref{eq32}, ce qui donne
\begin{equation} \label{eq36}
\eps_1 \epssoul^{\ssoul+1}  \frac{1}{(\jprun +1)^{s_1}}  \sum_{N \geq k \geq \ell  \geq 1} \Big( \frac{1}{(\ell+j'_2)^{s_2}(k+j'_3)^{s_3}} - \frac{1}{(k+j'_2)^{s_2}(\ell+j'_3)^{s_3}} \Big).
\end{equation} 
\item On regroupe   le deuxi\`eme terme de \eqref{eq31} avec le premier terme  de \eqref{eq33} ; en 
d\'ecouplant la sommation sur $k$ et $\ell$, on obtient en omettant le terme d'erreur 
\begin{multline}
\eps_1 \epssoul^{\ssoul+1} 
  \sum_{ k=1} ^N   \sum_{ \ell=1} ^N  \frac{1}{(k+\jprun)^{s_1}(\ell+j'_2)^{s_2}(k+j'_3)^{s_3}}  
\\
+ \eps_1 \epssoul^{\ssoul+1} 
  \sum_{ k=1} ^N  \frac{1}{(k+\jprun)^{s_1}(k+j'_2)^{s_2}(k+j'_3)^{s_3}}.  \label{eq37}
\end{multline} 
\item On regroupe   le deuxi\`eme terme de \eqref{eq35} avec le premier terme  
de \eqref{eq34} ; en d\'ecouplant la sommation, on obtient (en omettant le terme d'erreur)
\begin{multline}
-\eps_1 \epssoul^{\ssoul+1} 
  \sum_{ k=1} ^N   \sum_{ \ell=1} ^N  \frac{1}{(k+\jprun)^{s_1}(k+j'_2)^{s_2}(\ell+j'_3)^{s_3}}  \label{eq38}\\
- \eps_1 \epssoul^{\ssoul+1} 
  \sum_{ k=1} ^N  \frac{1}{(k+\jprun)^{s_1}(k+j'_2)^{s_2}(k+j'_3)^{s_3}}.
\end{multline} 
\item On regroupe   le premier  terme de \eqref{eq31} avec celui de \eqref{eq32}, d'o\`u :
\begin{multline} 
-\eps_1 \epssoul^{\ssoul+1} 
  \sum_{ k=1} ^N   \sum_{ \ell=1} ^N  \frac{1}{(k+\jprun+1)^{s_1}(k+j'_2)^{s_2}(\ell+j'_3)^{s_3}} \label{eq39} \\
- \eps_1 \epssoul^{\ssoul+1} 
  \sum_{ k=1} ^N  \frac{1}{(k+\jprun+1)^{s_1}(k+j'_2)^{s_2}(k+j'_3)^{s_3}}, 
\end{multline} 
qui se trouve \^etre la m\^eme \'equation que \eqref{eq38} mais avec $j'_1$ remplac\'e par $j'_1+1$ (et sans terme d'erreur \`a omettre).
\item On regroupe   le premier  terme de \eqref{eq30} avec celui de \eqref{eq35}, d'o\`u :
\begin{multline}
\eps_1 \epssoul^{\ssoul+1} 
  \sum_{ k=1} ^N   \sum_{ \ell=1} ^N  \frac{1}{(k+\jprun+1)^{s_1}(\ell+j'_2)^{s_2}(k+j'_3)^{s_3}}  \label{eq40}\\
+ \eps_1 \epssoul^{\ssoul+1}
  \sum_{ k=1} ^N  \frac{1}{(k+\jprun+1)^{s_1}(k+j'_2)^{s_2}(k+j'_3)^{s_3}} 
\end{multline} 
qui est la m\^eme \'equation que \eqref{eq37} mais avec $j'_1$ remplac\'e par $j'_1+1$ (et sans terme d'erreur \`a omettre).
\end{enumerate} 

\bigskip

Pour parvenir \`a la conclusion cherch\'ee, il suffit d'effectuer les groupements suivants, et de constater que chacun d'eux est de la forme voulue : 
\begin{itemize}
\item Pour tout $\eps_1 \in \zdz$, on regroupe la somme \eqref{eq36} correspondant aux triplets $(\eps_1, 1, 1)$, $(\eps_1, 1, -1)$, $(\eps_1, -1, 1)$ et  $(\eps_1, -1, -1)$. La somme de ces quatre termes vaut 
\begin{multline*} 
\frac{\eps_1^{s_1} }{(\jprun +1)^{s_1}}   
\bigg(\sum_{\eps_2, \eps_3 \in \zdz} \eps_2^{s_2+1} \eps_3^{s_3+1}
\\
\cdot \sum_{N \geq k \geq \ell  \geq 1} \Big( \frac{1}{(\ell+\eps_2 
\cdot j_2)^{s_2}(k+\eps_3 \cdot j_3)^{s_3}} - \frac{1}{(k+\eps_2 \cdot j_2)^{s_2}(\ell+\eps_3 \cdot j_3)^{s_3}} \Big)\bigg)
\end{multline*}
Le th\'eor\`eme \ref{thconj3} (d\'emontr\'e en profondeur 2 au \S \ref{subsec52}) s'applique \`a cette somme, et montre qu'elle s'\'ecrit $Q (H_N) + \gdo_\eps (N^{-1+\eps})$, o\`u $Q(0)$ est une combinaison lin\'eaire (\`a coefficients dans $\dd_n ^{-(s_2+s_3)} \Z$) de 1, de valeurs de $\zeta$ en des entiers impairs $s$ compris entre 3 et $s_2+s_3$, et de $\zetaetoile(s_3, s_2) - \zetaetoile(s_2, s_3)$. En outre ce polyz\^eta antisym\'etrique appara\^{\i}t avec un coefficient nul si $s_2$ ou $s_3$ est pair, et avec un coefficient
$4 \eps_1^{s_1} \frac{1}{(\jprun +1)^{s_1}}   $ si $s_2$ et $s_3$ sont impairs. Dans ce dernier cas, en sommant sur $\eps_1 \in \zdz$ on obtient finalement un coefficient
$$4 \Big(  \frac{1}{(j_1 +1)^{s_1}} +   \frac{(-1)^{s_1}}{(n-j_1) ^{s_1}}  \Big)$$
qui permet de justifier la remarque qui suit l'\'enonc\'e du th\'eor\`eme.
\item Pour tout $(\eps_1, \eps_3)  \in \zdzd$, on regroupe la somme double de \eqref{eq37} pour $(\eps_1, 1 , \eps_3)$ avec celle pour $(\eps_1, -1 , \eps_3)$, et avec  la somme double de \eqref{eq40} relative \`a $(-\eps_1, 1 , -\eps_3)$ et celle relative \`a $(-\eps_1,-1, -\eps_3)$. La contribution globale de ces 4 sommes doubles est, en notant g\'en\'eriquement $(\eta_1\eps_1, \eta_2, \eta_1\eps_3)$ les quatre triplets $\epssoul$ qui interviennent   :  
\begin{multline*}
\eps_1^{s_1}  \eps_3^{s_3+1} \bigg(
\sum_{\eta_1, \eta_2 \in \zdz} 
\eta_1^{s_1+s_3+1} \eta_2 ^{s_2+1} \\
\cdot 
  \sum_{ k=1} ^N   \sum_{ \ell=1} ^N  \frac{1}{(k+\eta_1 \cdot (\jprun))^{s_1}(\ell+\eta_2 \cdot j_2)^{s_2}(k+\eta_1 \cdot j'_3)^{s_3}}\bigg). 
\end{multline*}
Cette somme double se scinde sous la forme suivante :
\begin{multline} \label{eq41}
\eps_1^{s_1}  \eps_3^{s_3+1} 
\Big(   \sum_{ k=1} ^N \sum_{\eta_1  \in \zdz} 
 \frac{\eta_1^{s_1+s_3+1}  }{(k+\eta_1 \cdot (\jprun))^{s_1} (k+\eta_1 \cdot j'_3)^{s_3}} \Big)
\\
\cdot \Big(    \sum_{ \ell=1} ^N \sum_{\eta_2 \in \zdz} 
   \frac{  \eta_2 ^{s_2+1}}{ (\ell+\eta_2 \cdot j_2)^{s_2} } \Big). 
\end{multline} 
D'apr\`es le th\'eor\`eme \ref{thconj3} (d\'emontr\'ee en profondeur 1), la deuxi\`eme somme s'\'ecrit sous la forme $A_1(H_N) + \gdoeps$ o\`u $A_1$ est un polyn\^ome tel que $A_1(0)$ soit une combinaison lin\'eaire de 1 et de valeurs de $\zeta$ en des entiers $s$ impairs compris entre 3 et $s_2$, puisque $\zetaetoile(1)=0$. En outre $\dd_n^{s_2}$ est un d\'enominateur commun des coefficients de cette combinaison lin\'eaire. Enfin on a d\'emontr\'e au paragraphe \ref{subsec51} que $A_1(0) \in \Q$ si $s_2$ est pair, et $A_1(0) \in \Q + \Q \zeta(s_2)$ si $s_2$ est impair ; mais cette pr\'ecision suppl\'ementaire est inutile ici.

Pour la premi\`ere somme de \eqref{eq41},  on applique le th\'eor\`eme \ref{thconj2}, d\'emontr\'e en profondeur 1 (voir \S \S  \ref{subsec41} et \ref{subsec51}) : cette somme s'\'ecrit  sous la forme $A_2(H_N) + \gdoeps$ o\`u $A_2$ est un polyn\^ome tel que $A_2(0)$ soit une combinaison lin\'eaire de 1 et de valeurs de $\zeta$ en des entiers $s$ impairs compris entre 3 et $s_1+s_3$. En outre $\dd_n ^{s_1+s_3}$ est un d\'enominateur commun des coefficients de cette combinaison lin\'eaire. 

Comme la divergence logarithmique de $H_N$ est compens\'ee par le $N^\eps$ du terme d'erreur, on peut faire le produit des deux expressions pr\'ec\'edentes et obtenir
$$ \eqref{eq41}  = \eps_1 \epssoul^{\ssoul+1}  (A_1A_2)(H_N) + \gdoeps . $$
En outre, $A_1A_2(0)$ est une combinaison lin\'eaire de termes de la forme 1, $\zeta(s')$,  $\zeta(s'')$ ou  $\zeta(s') \zeta(s'')$, avec $s'$, $s''$ impairs et $3 \leq s' \leq s_1+s_3$, $3 \leq s'' \leq s_2$ ; et $\dd_n ^{s_1+s_2+s_3}$ est un d\'enominateur commun des coefficients. 
\item  Pour tout $(\eps_1, \eps_2)  \in \zdzd$, on regroupe la somme double de \eqref{eq38} pour $(\eps_1,  \eps_2, 1)$ avec celle pour $(\eps_1, \eps_2, -1)$, et avec  la somme double de \eqref{eq39} relative \`a $(-\eps_1 , -\eps_2, 1)$ et celle relative \`a $(-\eps_1, -\eps_2, -1)$. Le m\^eme ph\'enom\`ene que pr\'ec\'edemment se produit.
\item Pour tout $\epssoul \in \zdzt$, la somme simple de \eqref{eq37} et celle de \eqref{eq38} (pour cette m\^eme valeur de $\epssoul$) sont oppos\'ees donc leurs contributions \`a \eqref{eq12nv} s'annulent.
\item De m\^eme, la somme simple de \eqref{eq39} et celle de \eqref{eq40} s'annulent pour tout $\epssoul \in \zdzt$.
 \end{itemize}

\subsection{Preuve du th\'eor\`eme \ref{thconj3}  en profondeur quelconque}  \label{subsec54}

Dans ce paragraphe, on d\'emontre le th\'eor\`eme \ref{thconj3} en profondeur $p \geq 4$ en supposant (par r\'ecurrence) qu'il est vrai en profondeurs $p-2$ et $p-1$. En fait cette preuve fonctionne aussi quand $p=2$ et $p=3$ ; on retrouve alors les d\'emonstrations des deux paragraphes pr\'ec\'edents, \`a condition d'\^etre attentif aux conventions quand on somme sur des ensembles vides.
Notamment, \`a la convention habituelle
$$\sum_{k \in \emptyset } f(k) = 0$$
on adjoint la convention 
$$\sum_{k_1 \geq \ldots \geq k_r \geq 1} f(k_1, \ldots, k_r) = 1 \mbox{ pour } r = 0$$
car cette somme porte sur un ensemble vide de variables (par opposition \`a la pr\'ec\'edente, o\`u une variable parcourait un ensemble vide). 

\bigskip

L'initialisation de la r\'ecurrence se fait de mani\`ere tout \`a fait analogue au cas des profondeurs 2 et 3 : puisque $n$ est suppos\'e pair (voir la remarque \ref{rqparite}), il suffit, apr\`es avoir pos\'e 
$$\tau_{s_1,\ldots,s_p} = \sum_{N \geq k_1 \geq \ldots  \geq k_p \geq 1} \frac{1}{(k_1 + \frac{n}{2})^{s_1}  \ldots   (k_p + \frac{n}{2})^{s_p}   }, $$
de constater que l'on a 
$$
\tau_{s_1, \ldots,  s_p} 
= \sum_{p' = 0} ^p  (-1)^{p'} \Big( \sum_{\frac{n}{2} \geq \ell_p >  \ldots > \ell_{p-p'+1} \geq 1} \frac{1}{\ell_p ^{s_p}  \ldots  \ell_{p-p'+1} ^{s_{p-p'+1}}} \Big)
 \Big( \sum_{N+\frac{n}{2} \geq \ell_1 \geq  \ldots \geq \ell_{p-p'} \geq 1} \frac{1}{\ell_1 ^{s_1}  \ldots  \ell_{p-p'} ^{s_{p-p'}}} \Big). 
$$

D\'emontrons maintenant l'h\'er\'edit\'e, qui est la partie difficile. On suppose pour cela que 
le th\'eor\`eme \ref{thconj3} est vrai en profondeurs $p-2$ et $p-1$. 
On   adopte les notations suivantes : $\jsoul = (j_1, \ldots , j_p)$, $\ssoul = (s_1, \ldots,  s_p)$,  $\epssoul = (\eps_1, \ldots , \eps_p)$, 
$\epssoul \cdot \jsoul = (\eps_1 \cdot j_1, \ldots , \eps_p \cdot j_p)$, 
$ \epssoul ^{\ssoul + 1}  = \eps_1 ^{s_1+ 1} \ldots  \eps_p ^{s_p+ 1} $, 
$$K_N(\jsoul, \ssoul) = \sum_{N \geq k_1 \geq \ldots  \geq k_p  \geq 1}  \frac{1}{(k_1+j_1)^{s_1} \ldots   (k_p+j_p)^{s_p}}$$
et, pour $\sigma \in \spp$ : 
$$K_N^\sigma (\jsoul, \ssoul)  = K_N (j_{\sigma(1)},\ldots , j_{\sigma(p)}, s_{\sigma(1)}, \ldots , s_{\sigma(p)})$$
de telle sorte que $K_N^{\Id} (\jsoul, \ssoul) = K_N  (\jsoul, \ssoul) $. Comme $\eps_{\sigma^{-1}} = \eps_{\sigma}$, on a : 
$$\eqref{eq4} = \sum_{\epssoul \in \zdzp} \epssoul ^{\ssoul + 1} \sum_{\sigma \in \spp} \eps_\sigma
K_N^\sigma (\epssoul \cdot \jsoul, \ssoul).$$
On pose aussi 
$$\Delta_{\epssoul} ^\sigma (\jsoul) = 
K_N^\sigma (\eps_1 \cdot (j_1+1), \eps_2 \cdot j_2, \ldots,  \eps_p \cdot j_p, \ssoul)
- K_N^\sigma (\epssoul  \cdot \jsoul , \ssoul).$$
Alors la diff\'erence entre \eqref{eq4} pour $(j_1+1, j_2, \ldots , j_p)$ et \eqref{eq4} pour $\jsoul $ est :
\begin{equation} \label{eq60} 
\sum_{\epssoul \in \zdzp} \epssoul ^{\ssoul + 1} \sum_{\sigma \in \spp} \eps_\sigma \Delta_{\epssoul} ^\sigma (\jsoul) .
\end{equation}
La suite de la preuve est consacr\'ee \`a \eqref{eq60} : il s'agit de montrer que cette somme est de la forme voulue, ce qui terminera la r\'ecurrence (de m\^eme qu'en profondeur 2 et 3). 

\bigskip

Pour tout $\sigma \in \spp$, on pose $\tsig= \sigma^{-1}(1)$ et $j'_1 = \eps_1 \cdot j_1$, \ldots,   $j'_p = \eps_p \cdot j_p$, de telle sorte que $\jprsoul = (j'_1,\ldots, j'_p) =\epssoul \cdot \jsoul$. 
On pose aussi, par convention, $k_0 = N$ et $k_{p+1} = 1$. 
Supposons d'abord que $\eps_1 =  + 1$ ; on a dans ce cas $\eps_1 \cdot (j_1+1) = j'_1 + 1$, et : 
\begin{eqnarray*}
\desj 
&=& K_N^\sigma (j'_1+1, j'_2, \ldots,  j'_p, \ssoul) - K_N^\sigma (\jprsoul , \ssoul) \\
&=& \sum_{N \geq k_1 \geq \ldots  \geq k_p  \geq 1} (k_1+j'_{\sigma(1)})^{-s_{\sigma(1)}} \ldots   
 (k_{\tsig}+j'_1+1)^{-s_1} \ldots  (k_p+j'_{\sigma(p)})^{-s_{\sigma(p)}} \\
&& \quad \quad  - \sum_{N \geq k_1 \geq \ldots  \geq k_p  \geq 1} (k_1+j'_{\sigma(1)})^{-s_{\sigma(1)}} \ldots   
 (k_{\tsig}+j'_1)^{-s_1} \ldots  (k_p+j'_{\sigma(p)})^{-s_{\sigma(p)}} \\
&=&  \sum_{N \geq k_1 \geq \ldots \geq \widehat{k_{\tsig}} \geq \ldots  \geq k_p  \geq 1} 
(k_1+j'_{\sigma(1)})^{-s_{\sigma(1)}} \ldots   \widehat{ (k_{\tsig}+j'_1)^{-s_1} } \ldots  (k_p+j'_{\sigma(p)})^{-s_{\sigma(p)}} \\
&& \quad \quad \times  \sum_{k_{\tsig} = k_{\tsig+1}} ^{k_{\tsig-1}}  (k_{\tsig}+j'_1+1)^{-s_1} 
- (k_{\tsig}+j'_1)^{-s_1} \\
&=&  \sum_{N \geq k_1 \geq \ldots \geq  \widehat{k_{\tsig}} \geq \ldots  \geq k_p  \geq 1} 
(k_1+j'_{\sigma(1)})^{-s_{\sigma(1)}} \ldots   \widehat{ (k_{\tsig}+j'_1)^{-s_1} } \ldots  (k_p+j'_{\sigma(p)})^{-s_{\sigma(p)}} \\
&& \quad \quad \times  \Big( \frac{1}{ (k_{\tsig-1}+j'_1+1)^{s_1}} -  \frac{1}{ (k_{\tsig+1}+j'_1)^{s_1}} . 
\Big)
\end{eqnarray*}
Dans ce calcul, comme dans toute la suite, on note avec un chapeau l'omission d'un terme dans une liste. En outre, on utilise les conventions $k_0 = N$ et $k_{p+1} = 1$. 

Dans le cas o\`u  $\eps_1=-1$, la derni\`ere formule obtenue pour $\desj$ reste valable, \`a condition 
d'en prendre l'oppos\'e et d'y remplacer $j'_1$ par $j'_1-1$. Cela montre qu'on peut \'ecrire, quelle que soit la valeur de $\eps_1$ : 
\begin{equation} \label{eq59}
\desj = \sesj - \sesjti
\end{equation}
en posant 
\begin{eqnarray} 
\lefteqn{\sesj  } \nonumber \\ 
&=&\eps_1  \sum_{N \geq k_1 \geq \ldots \geq \widehat{k_{\tsig}} \geq \ldots  \geq k_p  \geq 1} 
(k_1+j'_{\sigma(1)})^{-s_{\sigma(1)}} \ldots   \widehat{ (k_{\tsig}+j'_1)^{-s_1} } \ldots  (k_p+j'_{\sigma(p)})^{-s_{\sigma(p)}}   \label{eq57} \\
&& \quad \quad \quad \quad \times    \frac{1}{ (k_{\tsig-1}+\jprun+1)^{s_1}}  \nonumber 
\end{eqnarray}
et 
\begin{eqnarray} 
\lefteqn{\sesjti } \nonumber \\
&=& \eps_1  \sum_{N \geq k_1 \geq \ldots \geq \widehat{k_{\tsig}} \geq \ldots  \geq k_p  \geq 1} 
(k_1+j'_{\sigma(1)})^{-s_{\sigma(1)}} \ldots   \widehat{ (k_{\tsig}+j'_1)^{-s_1} } \ldots  (k_p+j'_{\sigma(p)})^{-s_{\sigma(p)}}   \label{eq58} \\
&& \quad \quad \quad \quad \times    \frac{1}{ (k_{\tsig+1}+\jprun)^{s_1}}.  \nonumber 
\end{eqnarray}
La relation \eqref{eq59} va nous permettre de d\'emontrer que \eqref{eq60} est de la forme voulue. Dans un premier temps, on isole deux cas particuliers. Le premier concerne les termes de la forme $\sesj$ correspondant \`a des permutations $\sigma$ telles que $\tsig = 1$. Pour ces termes, on a d'apr\`es  \eqref{eq57} la majoration $\sesj = \gdo(\frac{(\log N)^{p-1}}{N})$ puisque $k_0 = N$ ; donc ces termes rentrent dans le terme d'erreur, et on peut les ignorer. Par ailleurs, si on regroupe tous les termes de la forme $\sesjti$ correspondant \`a des permutations $\sigma$ telles que $\tsig = p$, on obtient pour contribution globale \`a \eqref{eq60}, puisque $k_{p+1}=1$  : 
\begin{eqnarray} 
&& \frac{- 1}{ (\jprun + 1)^{s_1}} \sum_{\epssoul \in \zdzp} \epssoul ^{\ssoul + 1} \eps_1  \sum_\indso{\sigma \in \spp}{\tsig = p} \eps_\sigma   \label{eq62} \\
&&\times \sum_{N \geq k_1 \geq \ldots \geq k_{p-1}  \geq 1} 
(k_1+j'_{\sigma(1)})^{-s_{\sigma(1)}} \ldots (k_{p-1}+j'_{\sigma(p-1)})^{-s_{\sigma(p-1)}}.  \nonumber 
\end{eqnarray}
En fixant $\eps_1$ dans cette somme, on peut appliquer le th\'eor\`eme \ref{thconj3} en profondeur $p-1$, avec $(\eps_2, \ldots, \eps_p)$, $(j_2, \ldots, j_p)$, et $(s_2, \ldots, s_p)$. Le terme obtenu est multipli\'e par le rationnel $ \frac{- \eps_1 ^{s_1}}{ (\jprun + 1)^{s_1}} $, dont $\dd_n^{s_1}$ est un d\'enominateur ; le r\'esultat est donc de la forme souhait\'ee. Ce raisonnement g\'en\'eralise celui qui a permis, en profondeur 3, de traiter la somme \eqref{eq36}. 

Pour terminer la preuve, on peut donc ignorer dans \eqref{eq60} les termes provenant de ces deux familles de cas particuliers. Cela revient \`a faire la convention suivante, que nous adoptons dans toute la suite : 
\begin{equation} \label{eq61}
\left\{
\begin{array}{l}
\sesj = 0 \mbox{ si } \tsig = 1 \\
\sesjti = 0 \mbox{ si } \tsig = p. 
\end{array}
\right.
\end{equation}

On peut maintenant relier les sommes $\sesj$ et $\sesjti$, pour les \'etudier simultan\'ement. Pour cela, on d\'emontre l'\'egalit\'e suivante, valable pour tout $\sigma \in \spp$ tel que $\tsig \geq 2$ : 
\begin{equation} \label{eq56}
\sesj = \sesjtirondjtisoul
\end{equation}
avec $\jtisoul = (j_1+\eps_1, j_2, \ldots, j_p)$. Posons $\sigmati = \sigma \circ \transpotsig$ ; 
on a $\sigmati(j) = \sigma(j)$ pour $j \not\in \{\tsig-1, \tsig\}$, $\sigmati(\tsig-1) = 1$ et $\sigmati(\tsig) = \sigma(\tsig-1)$. En particulier, on a $\tsigti = \tsig - 1$. On constate alors qu'en rempla\c cant $j_1$ par $j_1 + \eps_1$ (ce qui revient \`a remplacer $j'_1$ par $j'_1+1$) dans la d\'efinition \eqref{eq58} de $\sesjtisigti$, on obtient exactement celle \eqref{eq57} de $\sesj$,  \`a un changement de notation pr\`es sur les indices de sommation. 
 En effet, dans \eqref{eq57}, l'indice $k_{\tsig}$ n'appara\^{\i}t pas dans la somme, alors que $k_{\tsig - 1}$ appara\^{\i}t et correspond \`a deux facteurs. Dans \eqref{eq58}, c'est $k_{\tsig-1}$ qui n'appara\^{\i}t pas, et $k_{\tsig}$ correspond \`a deux facteurs, qui sont exactement ceux provenant de $k_{\tsig-1}$ dans \eqref{eq57} (apr\`es avoir remplac\'e $j'_1$ par $j'_1+\eps_1$ dans \eqref{eq58}). Enfin les $k_j$ pour $j \not\in \{\tsig-1, \tsig\}$ jouent le m\^eme r\^ole dans \eqref{eq57} et dans \eqref{eq58}. Ceci termine la preuve de \eqref{eq56}. 

\bigskip

Compte tenu de \eqref{eq59}, \eqref{eq56} et \eqref{eq61}, on peut maintenant r\'e\'ecrire \eqref{eq60} sous la forme : 
\begin{equation} \label{eq63}
-  \sum_{\epssoul \in \zdzp} \epssoul ^{\ssoul + 1}   \sum_\indso{\sigma \in \spp}{\tsig \leq  p-1} \eps_\sigma  
\Big( \sesjti + \sesjtijti \Big)
\end{equation}
en omettant \eqref{eq62} et le terme d'erreur $ \gdo(\frac{(\log N)^{p-1}}{N})$ rencontr\'es plus haut (ce qui correspond \`a la convention \eqref{eq61}). 
Pour conclure la preuve, il suffit donc de d\'emontrer que \eqref{eq63} est de la forme voulue. 

\bigskip

Pour cela, on d\'efinit une application
\begin{eqnarray*}
\Phi : \{ \sigma \in \spp, \, \tsig \leq p-1 \} 
&\rightarrow & \unpmu \times \deuxp  \times \sppmd \\
\sigma &\mapsto& (\tsig, \tetasig, \gamma)
\end{eqnarray*}
de la fa\c con suivante. Pour $\sigma \in \spp$ tel que $\tsig \leq p-1$, on pose
$$\tetasig = \sigma(\tsig+1), $$
et on note $\phisig : \unpmd \rightarrow \unp \moins \{ \tsig, \tsig+1\}$ et 
$\psisig : \unpmd \rightarrow \deuxp \moins \{ \tetasig\}$ les bijections strictement croissantes. On pose alors
$$\gamma = \psisig^{-1} \circ \sigma \circ \phisig \in \sppmd = {\mathfrak S} (\unpmd)$$
o\`u on identifie $\sigma$ avec sa restriction $\sigma : \unp \moins \{ \tsig, \tsig+1\} \rightarrow \deuxp \moins \{ \tetasig\}$. Par d\'efinition de $\tsig$ et $\tetasig$, cette restriction est bijective, donc $\gamma$ aussi. Il est facile de voir que $\Phi$ est une bijection. 

\bigskip

Gr\^ace \`a cette bijection $\Phi$, on va remplacer la somme sur $\sigma $ dans \eqref{eq63} par une somme sur $(\tsig, \tetasig, \gamma)$. Pour cela on utilise la relation suivante, valable pour tout   $\sigma \in \spp$ tel que $\tsig \leq p-1$ : 
\begin{equation} \label{eq64}
\eps_\sigma  = (-1)^{\tetasig} \eps_\gamma.
\end{equation}
Pour d\'emontrer \eqref{eq64}, on \'etudie les couples $(i,j)$ tels que $1 \leq i <  j \leq p$ et $\sigma(i) > \sigma(j)$ ; la signature de $\sigma$ est donn\'ee par la parit\'e du nombre de tels couples. Soit $(i,j)$ un tel couple. Si $\{i,j\} \cap  \{ \tsig, \tsig+1\} = \emptyset$, ce couple correspond au couple $(\phisig^{-1}(i), \phisig^{-1}(j))$ qui contribue \`a la signature de $\gamma$. R\'eciproquement, chaque couple qui intervient dans le calcul de $\eps_\gamma$ est obtenu, une et une seule fois, de cette mani\`ere.  Comme le cas $\{i,j\} =  \{ \tsig, \tsig+1\}$
est exclu puisque $\sigma(\tsig) =1 < \sigma(\tsig+1)$, il y a exactement quatre autres possibilit\'es (qui s'excluent mutuellement) pour les couples $(i,j)$ qui contribuent \`a $\eps_\sigma$ mais pas \`a $\eps_\gamma$ : 
\begin{itemize}
\item Ou bien $i = \tsig$, mais c'est impossible car $\sigma(\tsig) = 1 < \sigma(j)$. 
\item Ou bien $i = \tsig  +1$ d'o\`u $j \geq \tsig  +  2$ avec $\sigma(j) < \tetasig$ ; le nombre de tels couples est $\card \{ j \geq \tsig + 2, \, \sigma(j) < \tetasig\}$.
\item Ou bien $j = \tsig$, d'o\`u $i < \tsig $ et $\sigma(i) > 1$ ; il y a exactement $\tsig - 1$ tels couples.
\item Ou bien $j = \tsig  +1$ d'o\`u $i < \tsig $ et $\sigma(i) > \tetasig$ ; le nombre de tels couples est $\card \{ i <   \tsig , \, \sigma(i) >  \tetasig\}$.
\end{itemize}
Pour d\'emontrer \eqref{eq64}, il suffit donc de prouver la relation suivante : 
\begin{equation} \label{eq65}
\card \{ j \geq \tsig + 2, \, \sigma(j) < \tetasig\} + \card \{ i <   \tsig , \, \sigma(i) >  \tetasig\}
+ \tsig - 1 \equiv \tetasig \mod 2. 
\end{equation}
Or on a clairement 
$$
\card \{ j \geq \tsig + 2, \, \sigma(j) < \tetasig\} + \card  \{ j \geq \tsig + 2, \, \sigma(j) >  \tetasig\} 
= \card\{ \tsig+2, \ldots, p \} = p - \tsig - 1 
$$
et
$$
\card  \{ i <   \tsig , \, \sigma(i) >  \tetasig\} + \card  \{ i \geq \tsig + 2, \, \sigma(i) >  \tetasig\} 
= \card\{ \tetasig+1, \ldots, p \} = p - \tetasig.   
$$
En additionnant ces deux relations on obtient \eqref{eq65}, ce qui termine la preuve de \eqref{eq64}.

\bigskip

Gr\^ace \`a la bijection $\Phi$ et \`a \eqref{eq64}, on peut maintenant \'ecrire  \eqref{eq63} sous la forme suivante : 
\begin{eqnarray}
\eqref{eq63}
&=& - \sum_{\teta = 2} ^p (-1)^\teta 
 \sum_{\epssoul \in \zdzp} \epssoul ^{\ssoul + 1} 
 \sum_{\gamma \in \sppmd} \eps_\gamma  \label{eq69} \\
&& \quad \quad \times \sum_{t = 1}^{p-1}  \Big( \setiphimu + \setiphimujti \Big) .  \nonumber 
\end{eqnarray}
On va maintenant montrer que la somme sur $t$ induit un d\'ecouplage de l'une des variables. Pr\'ecis\'ement, fixons $\teta \in \deuxp$, $\epssoul \in \zdzp$ et $\gamma \in \sppmd$. En posant $\sigt = \Phi^{-1}(t,\teta,\gamma)$ on a d'apr\`es \eqref{eq58} :
\begin{eqnarray*}
\sum_{t=1}^{p-1} \sesigtj
&=& \eps_1  \sum_{t=1}^{p-1}  \sum_{N \geq k_1 \geq \ldots \geq \widehat{k_t} \geq \ldots  \geq k_p  \geq 1} 
(k_{t+1}+\jprun)^{-s_1} 
\prod_\indso{1 \leq i \leq p }{i \neq t}  
(k_i+j'_{\sigt(i)})^{-s_{\sigt(i)}}.  
\end{eqnarray*}
Notons $\lam$ la variable $k_{t+1}$, qui appara\^{\i}t dans deux facteurs. Posons aussi $\ell_i = k_{\phisigt(i)}$ pour tout $i \in \unpmd$. On obtient : 
\begin{eqnarray*}
\sum_{t=1}^{p-1} \sesigtj
&=& \eps_1  \sum_{t=1}^{p-1}  \sum_{N \geq \ell_1 \geq \ldots \geq \ell_{t-1} \geq \lam \geq \ell_t \geq  \ldots  \geq \ell_{p-2}  \geq 1} 
(\lam +\jprun)^{-s_1}  (\lam + j'_{\tetasigt})^{-s_{\tetasigt}}  \\
&& \quad \quad \times 
\prod_{i=1} ^{p-2} (\ell_i+j'_{\sigt \circ \phisigt (i)})^{-s_{\sigt  \circ \phisigt (i)}}.  
\end{eqnarray*}
La propri\'et\'e cruciale est alors que le sommande est ind\'ependant de $t$, puisque $\tetasigt = \teta$ et $\sigt \circ  \phisigt = \psisigt \circ \gamma$  par d\'efinition  ; en outre $\psisigt$ ne d\'epend pas de $t$, mais seulement de $\teta$ (on note d\'esormais $\psii$ cette fonction). 
On peut donc d\'ecoupler la somme en \'ecrivant : 
$$\sum_{t=1}^{p-1}      
\sum_{N \geq \ell_1 \geq \ldots \geq \ell_{t-1} \geq \lam \geq \ell_t \geq  \ldots  \geq \ell_{p-2}  \geq 1} 
=
\sum_\indso{N \geq \ell_1 \geq \ldots  \geq \ell_{p-2}  \geq 1}{N \geq \lam \geq 1 }
+  \sum_{i=1}^{p-2}      \sum_\indso{N \geq \ell_1 \geq \ldots  \geq \ell_{p-2}  \geq 1}{\lam  = \ell_i }. 
$$
On obtient ainsi 
\begin{equation} \label{eq66}
\sum_{t=1}^{p-1} \sesigtj = \astgj +   \sum_{i=1}^{p-2}      \bstgj
\end{equation}
en posant 
\begin{eqnarray} 
\astgj & = &  \eps_1 
\Big( \sum_{\lam = 1}^N (\lam +\jprun)^{-s_1}  (\lam + j'_{\teta})^{-s_{\teta}}  \Big) \label{eq67} \\
&& \times \Big(    \sum_{N \geq \ell_1 \geq \ldots  \geq \ell_{p-2}  \geq 1}
\prod_{i=1} ^{p-2} (\ell_i+j'_{\psii \circ \gamma (i)})^{-s_{\psii \circ \gamma (i)}} \Big)     \nonumber
\end{eqnarray} 
et
\begin{eqnarray} 
\lefteqn{\bstgj}  \nonumber 
\\&=&   \eps_1   \sum_{N \geq \ell_1 \geq \ldots  \geq \ell_{p-2}  \geq 1}
(\ell_i +\jprun)^{-s_1}  (\ell_i + j'_{\teta})^{-s_{\teta}}   (\ell_i+j'_{\psii \circ \gamma (i)})^{-s_{\psii \circ \gamma (i)}}   \label{eq68} \\
&& \quad \quad \times   \prod_\indso{1 \leq i' \leq p-2}{i' \neq i}      (\ell_{i'}+j'_{\psii \circ \gamma (i')})^{-s_{\psii \circ \gamma (i')}}  .   \nonumber
\end{eqnarray} 
Gr\^ace \`a \eqref{eq66}, on peut maintenant \'ecrire \eqref{eq63} sous la forme suivante (en rempla\c cant dans \eqref{eq69}) : 
\begin{multline}
\eqref{eq63} = - \sum_{\teta = 2} ^p (-1)^\teta 
 \sum_{\epssoul \in \zdzp} \epssoul ^{\ssoul + 1} 
 \sum_{\gamma \in \sppmd} \eps_\gamma  \label{eq70} \\
  \times   \Big( \astgj +  \astgjti +    \sum_{i=1}^{p-2}    \Big(  \bstgj + \bstgjti \Big)  \Big). 
\end{multline}
Ici, les termes $\astgseul$ correspondent (en profondeur $p= 3$) aux sommes doubles des \'equations \eqref{eq37} \`a \eqref{eq40} ; les termes $\bstgseul$ correspondent aux sommes simples qui les accompagnent. 
On va maintenant g\'en\'eraliser le groupement de termes utilis\'e en profondeur 3 : ainsi, on groupe  les termes de \eqref{eq70} de telle sorte que chaque groupe soit de la forme voulue. Cela terminera la preuve du th\'eor\`eme \ref{thconj3}. 

\bigskip

La premi\`ere famille de groupements permet de traiter les termes $\astgseul$. Soient $\teta \in \deuxp$ et $(\epsunz, \epstetaz) \in \zdzd$ fix\'es. On regroupe les $2^{p-1}(p-2)!$ termes suivants :
\begin{eqnarray*}
&& \astgj  \mbox{ pour } \gamma \in \sppmd \mbox{ et  $\epssoul$ de la forme }
(\epsunz, \eta_2, \ldots, \eta_{\teta-1}, \epstetaz, \eta_{\teta+1}, \ldots, \eta_p) \\
&& \quad \quad \quad \mbox{ avec } (\eta_2, \ldots, \eta_{\teta-1},  \eta_{\teta+1}, \ldots, \eta_p) \in \zdzpmd \mbox{ , et} \\
&&\astgjti  \mbox{ pour } \gamma \in \sppmd \mbox{ et  $\epssoul$ de la forme }
(-\epsunz, \eta_2, \ldots, \eta_{\teta-1}, -\epstetaz, \eta_{\teta+1}, \ldots, \eta_p) \\
&& \quad \quad \quad \mbox{ avec } (\eta_2, \ldots, \eta_{\teta-1},  \eta_{\teta+1}, \ldots, \eta_p) \in \zdzpmd .
\end{eqnarray*}
Pour unifier ces deux cas, on note $\eps_1 = \eta_1 \epsunz$ et $\eps_\teta = \eta_1 \epstetaz$ avec $\eta_1 \in \zdz$. 
Pour les $2^{p-2}(p-2)!$ termes qui correspondent \`a $\eta_1 = -1$ (c'est-\`a-dire ceux de la forme $\astgjti$), il convient de remarquer qu'on a $\eps_1 = - \epsunz$ donc $\eps_1 \cdot (j_1+\eps_1) + \frac{\eps_1-1}{2} = (-1) \cdot ((\epsunz \cdot j_1) + \frac{\epsunz-1}{2})$. Ceci permet de prouver  que la contribution globale de ces $2^{p-1}(p-2)!$ termes \`a \eqref{eq70} s'\'ecrit, \`a un signe pr\`es qui d\'epend de $\teta$, $\epsunz$ et $\epstetaz$ : 
\begin{multline}
\Big( \sum_{\eta_1 \in \zdz} \eta_1 ^{s_1+s_\teta+1} \sum_{\lam=1}^N 
(\lam + \eta_1 \cdot ((\epsunz \cdot j_1)+  \frac{\epsunz-1}{2}))^{-s_1}
(\lam + \eta_1 \cdot (\epstetaz \cdot j_\teta))^{-s_\teta} \Big)  \label{eq67bis} 
\\
\times \Big( \sum_{(\eta_2, \ldots, \widehat{\eta_{\teta}}, \ldots, \eta_p) \in \zdzpmd}
\eta_2 ^{s_2+1} \ldots \widehat{\eta_\teta^{s_\teta+1}} \ldots \eta_p ^{s_p+1} 
 \sum_{\gamma \in \sppmd} \eps_\gamma   
\\
 \quad \quad \quad   \quad \quad   \times  
  \sum_{N \geq \ell_1 \geq \ldots  \geq \ell_{p-2}  \geq 1}
\prod_{i=1} ^{p-2} (\ell_i+\eta_{\psii \circ \gamma (i)} \cdot j_{\psii \circ \gamma (i)})^{-s_{\psii \circ \gamma (i)}} \Big).  
\end{multline}
Pour traiter le deuxi\`eme facteur de ce produit, on applique le th\'eor\`eme \ref{thconj3} en profondeur $p-2$, avec $j_2, \ldots, \widehat{j_\teta}, \ldots, j_p$ et $s_2, \ldots, \widehat{s_\teta}, \ldots, s_p$. Ce facteur s'\'ecrit donc  $A_1(H_N) + \gdoeps$ o\`u $A_1$ est un polyn\^ome tel que $A_1(0)$ soit une combinaison lin\'eaire de produits de la forme \eqref{eq999} avec $\{ i_1 , \ldots , i_{q-q'}  \} \cup \{ j_1 , \ldots , j_{2q'} \} \subset \deuxp \moins \{\teta\}$. De plus $\dd_n^{s_2 +  \ldots + \widehat{s_\teta} + \ldots + s_p}$ est un d\'enominateur commun des coefficients de cette combinaison lin\'eaire. 

Pour le premier facteur de \eqref{eq67bis},  on applique le th\'eor\`eme \ref{thconj2}, d\'emontr\'e en profondeur 1 (voir \S \S  \ref{subsec41} et \ref{subsec51}). Cette somme s'\'ecrit  donc sous la forme $A_2(H_N) + \gdoeps$ o\`u $A_2$ est un polyn\^ome tel que $A_2(0)$ soit une combinaison lin\'eaire de 1 et de valeurs de $\zeta$ en des entiers $s$ impairs compris entre 3 et $s_1+s_\teta$. En outre $\dd_n ^{s_1+s_\teta}$ est un d\'enominateur commun des coefficients de cette combinaison lin\'eaire. 

Comme la divergence logarithmique de $H_N$ est compens\'ee par le $N^\eps$ du terme d'erreur, on peut faire le produit des deux expressions pr\'ec\'edentes et obtenir
$$ \eqref{eq67bis}  =   (A_1A_2)(H_N) + \gdoeps . $$
En outre, $A_1A_2(0)$ est bien de la forme voulue. Ceci termine le traitement des termes de la forme  $\astgseul$ dans  \eqref{eq70}, car ces $2^{p+1}(p-1)!$ termes sont r\'epartis en $4(p-1)$ tels groupes. 

\bigskip

On va maintenant traiter les termes $\bstgseul$ de  \eqref{eq70}.
Pour cela, on les groupe deux par deux de la mani\`ere suivante. Soient $\teta \in \deuxp$,  $\epssoul \in \zdzp$, $\gamma \in \sppmd$ et $i \in \unpmd$ fix\'es. On note $\psiteta$ la bijection strictement croissante de $\unpmd$ dans $\deuxp \moins \{\teta\}$. Posons $\teta' = \psiteta(\gamma(i))$, $\alpha = \psitetapr^{-1}(\teta)$ et $\beta = \gamma(i) = \psiteta^{-1}(\teta')$. On note $\albe$ le cycle $\albelongun$ si $\alpha \leq \beta$, et le cycle $\albelongde$ si $\alpha > \beta$. On pose 
$\gamma' = \albe \circ \gamma$. Avec ces notations, on a $\gamma'(i) = \alpha$ d'o\`u $\{\teta, \psiteta(\gamma(i))\} = \{\teta' , \psitetapr(\gamma'(i))\}$. En outre, la d\'efinition de $\gamma'$ montre que pour tout $i' \in \unpmd \moins \{i\}$ on a $\psiteta(\gamma(i')) = \psitetapr(\gamma'(i'))$.  En reportant dans \eqref{eq68} on en d\'eduit : 
\begin{equation} \label{eq68bis}
\bstgj = \bstgjprime.
\end{equation}
Or on voit facilement que $\eps_{\gamma'} = \eps_{\gamma} (-1)^{\beta - \alpha} =  \eps_{\gamma} (-1)^{\teta - \teta'  - 1}$, d'o\`u $(-1)^\teta \eps_\gamma = - (-1)^{\teta'} \eps_{\gamma'}$. Donc les deux membres de l'\'egalit\'e \eqref{eq68bis} apparaissent dans \eqref{eq70} avec des signes oppos\'es : leurs contributions se neutralisent. Comme l'application $(\teta, \gamma) \mapsto (\teta', \gamma')$ ainsi d\'efinie est involutive, elle permet de grouper deux par deux tous les termes $\bstgj$ et $\bstgjti$ apparaissant dans \eqref{eq70}. Ceci d\'emontre que leur contribution globale est nulle, et termine la preuve du th\'eor\`eme \ref{thconj3}.

\section{Preuve du  th\'eor\`eme d\'ecoupl\'e}  \label{sec6}
 
D\'emontrons maintenant le th\'eor\`eme \ref{thdecouple}. La strat\'egie g\'en\'erale est la m\^eme que pour le th\'eor\`eme \ref{thconj1}, mais elle est beaucoup plus facile \`a mettre en \oe uvre.

\bigskip

Soit $P(k_1,\ldots,k_p)$ un polyn\^ome de degr\'e $\leq A(n+1)-2$ par rapport \`a chacune des variables. Comme au paragraphe \ref{subsec41}, on consid\`ere la  fraction rationnelle 
\begin{equation} \label{eq641}
R(k_1,\ldots, k_p) = \frac{P(k_1,\ldots, k_p)}{(k_1)_{n+1}^A \ldots (k_p)_{n+1}^A}
\end{equation}
dont la d\'ecomposition en \'el\'ements simples s'\'ecrit
\begin{equation} \label{eq662}
R(k_1,\ldots, k_p) = \sum_\indso{0 \leq j_1, \ldots, j_p \leq n}{1 \leq s_1, \ldots, s_p \leq A} 
\frac{\cjs}{(k_1+j_1)^{s_1} \ldots (k_p+j_p)^{s_p}}
\end{equation}
avec des rationnels $\cjs$. L'hypoth\`ese faite sur $P$ dans le  th\'eor\`eme \ref{thdecouple} s'\'ecrit
$$R(k_1, \ldots, k_{\ell-1}, -k_\ell-n, k_{\ell+1}, \ldots, k_p) = - R(k_1, \ldots, k_p)
\mbox{ pour tout } \ell \in \unp.$$
Par unicit\'e du d\'eveloppement en \'el\'ements simples, elle implique
\begin{equation} \label{eq287}
\cjsparsix = (-1)^{s_\ell+1} \cjs
\end{equation}
pour tous $j_1, \ldots, j_p, s_1, \ldots, s_p$ et pour tout $\ell \in \unp$. 

\bigskip

La s\'erie  \eqref{eqdecouple} est la limite, quand $N$ tend vers l'infini, de la somme 
\begin{equation} \label{eq288}
\sum_{k_1 = 1} ^N \ldots \sum_{k_p = 1} ^N R(k_1, \ldots, k_p).
\end{equation}
Pour tout entier $s \geq 1$, posons 
$$\zeta_N (s) = \sum_{k=1} ^N \frac{1}{k^s}.$$
Pour $s= 1$ c'est la somme harmonique (not\'ee aussi $H_N$), et pour $s \geq 2$ la suite $(\zeta_N (s))$ tend vers $\zeta(s)$ quand $N$ tend vers l'infini. On a, pour tous $(j_1, \ldots, j_p)$ et $(s_1, \ldots, s_p)$ : 
$$\sum_{1 \leq k_1, \ldots, k_p \leq N}\prod_{i=1} ^p   \frac{1}{(k_i + j_i)^{s_i}} 
= \prod_{i=1} ^p \Big( \zeta_{N+j_i}(s_i) - \sum_{k_i = 1} ^{j_i} \frac{1}{k_i^{s_i}} \Big).$$

Donc la somme \eqref{eq288} s'\'ecrit 
\begin{equation} \label{eq286}
  \sum_\indso{0 \leq j_1, \ldots, j_p \leq n}{1 \leq s_1, \ldots, s_p \leq A} 
 \cjs \prod_{i=1} ^p \Big( \zeta_{N+j_i}(s_i) - \sum_{k_i = 1} ^{j_i} \frac{1}{k_i^{s_i}} \Big).
  \end{equation}

Notons $E = \zeron^p \croix \una^p$ et consid\'erons la relation d'\'equivalence $\calr$ sur $E$ d\'efinie par :
$$
\begin{array}{l}
(j_1,\ldots, j_p, s_1, \ldots, s_p) \equiv (j'_1,\ldots, j'_p, s'_1, \ldots, s'_p) \, \bmod \, \calr \\
\\
\quad \quad \mbox{ si, et seulement si, } \\
\\
	\left\{
	\begin{array}{l}
	s_1 = s'_1, \ldots, s_p = s'_p\\
	j_1 \in \{j'_1, n-j'_1\}, \ldots, j_p \in \{j'_p, n-j'_p\}.
	\end{array}
	\right.
\end{array}
$$
On peut scinder la somme \eqref{eq286} en somme sur les classes d'\'equivalence\footnote{Il s'agit des orbites sous l'action de $\zdzp$ sur $E$ d\'efinie au paragraphe \ref{subsec41}.} modulo $\calr$ (puisque celles-ci forment une partition de $E$). Nous allons d\'emontrer que la somme sur chaque classe est de la forme $Q(H_N) + \odu$ o\`u $Q$ est un polyn\^ome, $H_N$ la somme harmonique et $\odu$ une suite qui tend vers 0, avec la propri\'et\'e que $Q(0)$ est un polyn\^ome  \`a coefficients rationnels, de degr\'e au plus $p$,  en les  $\zeta(s)$, pour $s$ entier impair compris entre 3 et $A$. Quand $N$ tend vers l'infini, la somme \eqref{eq286} converge vers \eqref{eqdecouple} donc la contribution globale de ces polyn\^omes $Q(H_N)$ sera un polyn\^ome constant, dont  la valeur (en 0) est de la forme annonc\'ee dans le  th\'eor\`eme \ref{thdecouple}. Ceci d\'emontrera donc  le  th\'eor\`eme \ref{thdecouple}.

\bigskip

D\'emontrons maintenant ce fait. 
Soit $(j_1,\ldots, j_p, s_1, \ldots, s_p) \in E$. Pour simplifier les notations, on suppose (quitte \`a permuter les indices) que $j_1  =  \ldots = j_a = \frac{n}{2}$ et que $j_{a+1}, \ldots, j_p$ sont diff\'erents de $n/2$, avec $ a \in \zerop$ (par exemple $a=0$ d\`es que $n$ est impair). Alors la classe d'\'equivalence de $(j_1,\ldots, j_p, s_1, \ldots, s_p)$ modulo $\calr$ est form\'ee par les $2^{p-a}$ \'el\'ements $(\frac{n}{2},\ldots, \frac{n}{2}, j'_{a+1}, \ldots, j'_p, s_1, \ldots, s_p)$ tels que $j'_{a+1} \in \{j_{a+1}, n-j_{a+1}\}$, \ldots, $ j'_p \in \{j_p, n-j_p\}$. Pour $\eps \in \{-1, 1\}$ et $j \in \zeron$ on pose (comme au paragraphe \ref{subsec41}) :
$$
\left\{
\begin{array}{l}
\eps \cdot j =j \mbox{ si } \eps = + 1,\\
\eps \cdot j =n-j \mbox{ si } \eps = - 1.
\end{array}
\right.
$$
Alors ces $2^{p-a}$ \'el\'ements s'\'ecrivent $(\eps_1 \cdot j_1, \ldots, \eps_p \cdot j_p, s_1, \ldots, s_p)$ o\`u $(\eps_1, \ldots, \eps_p)$ d\'ecrit $\{1\}^a \croix \{-1,1\}^{p-a}$ (c'est-\`a-dire que $\eps_1, \ldots, \eps_a$ valent toujours 1 et que $\eps_{a+1}, \ldots, \eps_p$ peuvent valoir 1 ou $-1$). La relation \eqref{eq287} donne alors, pour tout $ (\eps_1, \ldots, \eps_p) \in \{1\}^a \croix \{-1,1\}^{p-a}$  : 
$$
\cjsparsixbis = \eps_{a+1}^{s_{a+1}+1} \ldots  \eps_p^{s_p+1} \cjs,
$$
donc la somme \eqref{eq286} restreinte \`a la classe d'\'equivalence de $(j_1,\ldots, j_p, s_1, \ldots, s_p)$  est le produit de $\cjs$ par : 
\begin{eqnarray*}
\lefteqn{\sum_{(\eps_1, \ldots, \eps_p) \in \{1\}^a \croix \{-1,1\}^{p-a}} \eps_{a+1}^{s_{a+1}+1} \ldots  \eps_p^{s_p+1}
 \prod_{i=1} ^p \Big( \zeta_{N + \eps_i \cdot j_i} (s_i) - \sum_{k_i = 1} ^{\eps_i \cdot j_i} \frac{1}{k_i^{s_i}} \Big)} 
\nonumber \\
& = & \Big(   \prod_{i=1} ^a ( \zeta_{N + \frac{n}{2}} (s_i) - \sum_{k_i = 1} ^{n/2} \frac{1}{k_i^{s_i}} ) \Big)
 \sum_{(\eps_{a+1}, \ldots, \eps_p) \in \{-1,1\}^{p-a}}  \prod_{i=a+1} ^p \eps_i^{s_i+1}  ( \zeta_{N + \eps_i \cdot j_i} (s_i) - \sum_{k_i = 1} ^{\eps_i \cdot j_i} \frac{1}{k_i^{s_i}} ) \\
& = & \Big(   \prod_{i=1} ^a ( \zeta_{N + \frac{n}{2}} (s_i) - \sum_{k_i = 1} ^{n/2} \frac{1}{k_i^{s_i}} ) \Big) 
\prod_{i=a+1} ^p \sum_{\eps_i \in \{-1,1\}} \Big(  \eps_i^{s_i+1} \zeta_{N + \eps_i \cdot j_i} (s_i) -  \eps_i^{s_i+1} \sum_{k_i = 1} ^{\eps_i \cdot j_i} \frac{1}{k_i^{s_i}} \Big) \\
& = & \Big(   \prod_{i=1} ^a ( \zeta_{N + \frac{n}{2}} (s_i) - \sum_{k_i = 1} ^{n/2} \frac{1}{k_i^{s_i}} ) \Big) 
\prod_{i=a+1} ^p \Big( \zeta_{N+j_i}(s_i) +(-1)^{s_i+1}  \zeta_{N+n-j_i} (s_i) \\
&& \hspace{7cm}    -  \sum_{k_i = 1} ^{ j_i} \frac{1}{k_i^{s_i}} - (-1)^{s_i+1}   \sum_{k_i = 1} ^{n- j_i} \frac{1}{k_i^{s_i}}  \Big)\\
& = & \Big(   \prod_{i=1} ^a ( \zeta_{N} (s_i) - \sum_{k_i = 1} ^{n/2} \frac{1}{k_i^{s_i}} + \gdo(\frac{1}{N}) ) \Big) 
\prod_{i=a+1} ^p \Big( (1+(-1)^{s_i+1}) \zeta_N(s_i) + \gdo(\frac{1}{N})  \\
&& \hspace{7cm}    -  \sum_{k_i = 1} ^{ j_i} \frac{1}{k_i^{s_i}} - (-1)^{s_i+1}   \sum_{k_i = 1} ^{n- j_i} \frac{1}{k_i^{s_i}}  \Big), 
 \end{eqnarray*}
puisque $\zeta_{N+1}(s) = \zeta_N(s) +\gdo(1/N)$. Ce produit est bien de la forme $Q(H_N) + \odu$, o\`u $Q$ est un polyn\^ome (\`a coefficients r\'eels) tel que 
$$Q(0) =  \Big(   \prod_{i=1} ^a ( \zetaetoile (s_i) - \sum_{k_i = 1} ^{n/2} \frac{1}{k_i^{s_i}}) \Big) 
\prod_{i=a+1} ^p \Big( (1+(-1)^{s_i+1}) \zetaetoile (s_i)   -  \sum_{k_i = 1} ^{ j_i} \frac{1}{k_i^{s_i}} - (-1)^{s_i+1}   \sum_{k_i = 1} ^{n- j_i} \frac{1}{k_i^{s_i}}  \Big)$$
avec $\zetaetoile(1) = 0$ et $\zetaetoile(s) = \zeta(s)$ pour $s \geq 2$ (comme au  paragraphe \ref{subsec31}). 

\smallskip

Si l'un au moins parmi $s_1$, \ldots, $s_a$ est pair, alors la relation \eqref{eq287} montre que le coefficient 
$\cjs$ est nul, donc la classe d'\'equivalence de $(j_1,\ldots, j_p, s_1, 
\ldots, s_p)$  ne contribue pas \`a la somme \eqref{eq286}. On peut donc supposer 
que $s_1$, \ldots, $s_a$ sont tous impairs. Or l'expression  ci-dessus de $Q(0)$ ne 
fait appara\^{\i}tre, parmi les $\zeta(s_{i})$ avec $i \in \{a+1 , \ldots, p\}$, que 
ceux tels que $s_i$ soit impair ; en outre ceux parmi $s_1$, \ldots, $s_p$ qui valent 
1 disparaissent car $\zetaetoile(1) = 0$. Donc la contribution de  la classe d'\'equivalence 
de $(j_1,\ldots, j_p, s_1, \ldots, s_p)$  \`a la somme \eqref{eq286} est bien de la forme 
$Q(H_N)  + \odu$, o\`u $Q(0)$ est   un polyn\^ome  \`a coefficients rationnels, de degr\'e 
au plus $p$,  en les  $\zeta(s)$, pour $s$ entier impair compris entre 3 et $A$. Comme 
remarqu\'e ci-dessus, cela termine la preuve du th\'eor\`eme \ref{thdecouple}. 

\newcommand{\url}{\texttt}
\providecommand{\bysame}{\leavevmode ---\ }
\providecommand{\og}{``}
\providecommand{\fg}{''}
\providecommand{\smfandname}{\&}
\providecommand{\smfedsname}{\'eds.}
\providecommand{\smfedname}{\'ed.}
\providecommand{\smfmastersthesisname}{M\'emoire}
\providecommand{\smfphdthesisname}{Th\`ese}

\bigskip 

J. Cresson, 
Laboratoire de Math\'ematiques appliqu\'ees de Pau,
B\^atiment I.P.R.A, Univer\-si\-t\'e de Pau et des Pays de l'Adour,
avenue de l'Universit\'e, BP 1155, 64013 Pau cedex, France.

\medskip

S. Fischler,
Univ. Paris-Sud, Laboratoire de  Math\'ematiques,  UMR CNRS 8628, 
B\^atiment 425,
91405 Orsay cedex, France.

\medskip

T. Rivoal,
Institut Fourier,
CNRS UMR 5582, Universit{\'e} Grenoble 1,
100 rue des Maths, BP~74,
38402 Saint-Martin d'H{\`e}res cedex,
France.

\end{document}